\documentclass[11pt,reqno]{amsart}

\usepackage{amsmath}
\usepackage{amsthm}
\usepackage{graphicx}
\usepackage{amsfonts}
\usepackage{amssymb}
\usepackage{hyperref}
\usepackage{bm}

\usepackage[margin=1.25in]{geometry}

\usepackage{enumerate}
\numberwithin{equation}{section}

\newtheorem{theorem}{Theorem}[section]
\newtheorem{lemma}[theorem]{Lemma}
\newtheorem{proposition}[theorem]{Proposition}
\newtheorem{corollary}[theorem]{Corollary}

\theoremstyle{definition}
\newtheorem{definition}[theorem]{Definition}
\newtheorem{remark}[theorem]{Remark}

\def\E{{\mathbb E}}

\def\R{{\mathbb R}}
\def\N{{\mathbb N}}
\def\PP{{\mathbb P}}
\def\FF{{\mathbb F}}
\def\P{{\mathcal P}}
\def\Z{{\mathbb Z}}
\def\A{{\mathcal A}}
\def\F{{\mathcal F}}
\def\tr{{\mathrm{Tr}}}
\def\Aut{{\mathrm{Aut}}}
\def\Var{{\mathrm{Var}}}
\def\Cov{{\mathrm{Cov}}}
\def\deg{\mathrm{deg}}
\newcommand{\Erdos}{Erd\H{o}s-R\'enyi}
\newcommand{\lap}{L}
\def\newW{{\Lambda}}
\def\newA{{\Xi}}
\def\varfunc{{V}}

\title[Stochastic games on large graphs]{A case study on stochastic games on large graphs in mean field and sparse regimes}
\author{Daniel Lacker}
\author{Agathe Soret}
\address{Department of Industrial Engineering \& Operations Research, Columbia University}
\email{daniel.lacker@columbia.edu}
\email{acs2298@columbia.edu}

\begin{document}

\begin{abstract}
   We study a class of linear-quadratic stochastic differential games in which each player interacts directly only with its nearest neighbors in a given graph. We find a semi-explicit Markovian equilibrium for any transitive graph, in terms of the empirical eigenvalue distribution of the graph's normalized Laplacian matrix. This facilitates large-population asymptotics for various graph sequences, with several sparse and dense examples discussed in detail. In particular, the mean field game is the correct limit only in the dense graph case, i.e., when the degrees diverge in a suitable sense. Even though equilibrium strategies are nonlocal, depending on the behavior of all players, we use a correlation decay estimate to prove a propagation of chaos result in both the dense and sparse regimes, with the sparse case owing to the large distances between typical vertices. Without assuming the graphs are transitive, we show also that the mean field game solution can be used to construct decentralized approximate equilibria on any sufficiently dense graph sequence. 
\end{abstract}

\maketitle

\setcounter{tocdepth}{2}

\tableofcontents

\section{Introduction}

Mean field game (MFG) theory has enjoyed rapid development and widespread application since its introduction over a decade and a half ago by \cite{huang2006large,lasry-lions}. It provides a systematic framework for studying a broad class of stochastic dynamic games with many interacting players, in terms of limiting models featuring a continuum of players which are often more tractable. 
There are by now various rigorous results justifying the MFG approximation.
On the one hand, the equilibria of $n$-player games can be shown to converge to the MFG limit under suitable assumptions.
On the other hand, a solution of the continuum model may be used to construct \emph{approximate} equilibria for the $n$-player model with the particularly desirable properties of being \emph{decentralized} and \emph{symmetric}, in the sense that each player applies an identical feedback control which ignores the states of all other players.
We refer to the recent book of \cite{carmona-delarue-book} for a thorough account of MFG theory and its many applications.

A key structural assumption of the MFG paradigm is that the players interact \emph{symmetrically}, i.e., through an empirical measure which weights each player equally. In many natural situations, however, players do not view each other as exchangeable and instead interact directly only with certain subsets of players to which they are \emph{connected}, e.g., via some form of a graph or network. This is the purview of the broad field of \emph{network games}, and we refer to \cite{jackson2010social} for a representative overview of mostly static models.

Our paper contributes to a very recent line of work bridging MFG theory and network games by studying $n$-player stochastic dynamic games in which interactions are governed by a graph $G_n$ on $n$ vertices.
(When $G_n$ is the complete graph we recover the traditional MFG setting.) Roughly speaking, the goal is to understand the robustness of the mean field approximation and, when it fails, a substitute. Somewhat more precisely, two central questions are:
\begin{enumerate}
\item For what kinds of graph sequences $\{G_n\}$ is the usual MFG approximation still valid? 
\item What is the right \emph{limit model} for a general sequence $\{G_n\}$, and how well does it approximate the corresponding $n$-player game?
\end{enumerate}
Little progress has been made so far toward a systematic understanding of these questions.
The recent paper of \cite{delarue2017mean} addresses (1) when $G_n=G(n,p)$ is the \Erdos \ graph on $n$ vertices with fixed edge probability $p \in (0,1)$, showing that the usual MFG limit is still valid.
More recently, \cite{fouque-recent} study a linear-quadratic model very similar to ours, but only considering directed path or cycle graphs.
In another direction, recent efforts on (2) have proposed continuum models based on \emph{graphons}, which describe limit objects for general dense graph sequences  (\cite{lovasz2012large}). See recent work on \emph{graphon games} for the static case (\cite{carmona2019stochastic,parise2019graphon}) or \emph{graphon MFGs} in the dynamic case (\cite{caines2018graphon,gao2020lqg,vasal2020sequential}).

The combination of \emph{network}, \emph{dynamic}, and \emph{game-theoretic} effects is essential for many recent models of large economic and financial systems, and several recent studies have attacked specific models combining some of these features; see \cite{capponi2019dynamic,carmona2013mean,
feinstein2019dynamic,nadtochiy2018mean} and references therein.
Even without game-theoretic (strategic) features, incorporating network effects into large-scale dynamic models already presents many mathematical challenges, which very recent work has begun to address; see \cite{bhamidi2019weakly,coppini2019law,medvedev2014nonlinear} and \cite{detering2019directed,lacker2019large,oliveira2019interacting} for studies of dense and sparse graph regimes, respectively. Notably, prior work studied \emph{dense} graph regimes, and most questions in the \emph{sparse} regime (roughly defined as finite limiting neighborhood sizes) remain open, as was highlighted in particular in the recent paper of \cite{caines2018graphon}.

The purpose of our article is to give  comprehensive answers to (1) and (2), in both dense and sparse regimes, in the setting of a specific yet rich linear-quadratic model, inspired by the systemic risk (flocking) model of \cite{carmona2013mean}. 
For a suitably \emph{dense} graph sequence $\{G_n\}$ (meaning roughly that the degrees diverge as $n\to\infty$), we show (in Theorem \ref{th:approxeq}) that the classical construction of MFG theory is still valid: The MFG equilibrium gives rise to a sequence of \emph{decentralized and symmetric} approximate Nash equilibria for the $n$-player games. The dense regime includes the complete graph, the \Erdos \ graph $G_n=G(n,p_n)$ with $np_n \to \infty$, and many others.
Our findings in the dense case conform to an increasingly well-understood principle of statistical physics, that (static) interacting particle systems (e.g., the Ising model) on sufficiently dense and regular graphs tend to behave like their mean field counterparts (e.g., the Curie-Weiss model); see \cite{basak2017universality} and references therein.

The case of sparse graphs is more delicate, and the MFG approximation is no longer valid.
Here we restrict our attention to \emph{(vertex) transitive} graphs, which intuitively ``look the same" from the perspective of any vertex (see Definition \ref{def:transitive}); transitive graphs have rich enough symmetry groups to make up for the lack of exchangeability. We compute the Markovian Nash equilibrium explicitly (in Theorem \ref{thm_equilibrium_n_player}), up to the solution of a one-dimensional ordinary differential equation (ODE) governed by the empirical eigenvalue distribution of the Laplacian of the graph (i.e., the rate matrix of the simple random walk). As a consequence, we show (in Theorem \ref{thm:convergence_X^Gn_k}) that for a given graph sequence $\{G_n\}$, the limiting law can be computed for a typical player's state process, under the assumption that the empirical eigenvalue distributions of the Laplacian matrices of the graphs converge weakly. 
We also discuss (in Section \ref{se:cooperative}) similar and much simpler results for the corresponding cooperative problem, which we can explicitly solve for any graph (not necessarily transitive).

The eigenvalue distribution of a graph Laplacian is reasonably tractable in many interesting cases.
The dense graph case is precisely the case where the eigenvalue distribution converges weakly to a point mass.
In the sparse case, the eigenvalue distribution converges to a non-degenerate limit, and we characterize the much different $n\to\infty$ behavior in terms of this limit. 
We do not have a complete answer to (2) in the sparse case, as it remains unclear how to identify the limiting dynamics \emph{intrinsically}, without relying on $n\to\infty$ limits of $n$-player models. In contrast, the MFG framework identifies an intrinsic continuum model, the solution of which agrees with the $n\to\infty$ limit of the $n$-player equilibria.
See Section \ref{se:intrinsic} for further discussion of this point.

A crucial challenge in the sparse setting is that equilibrium controls are not \emph{local}, even at the limit. Even though each player's cost function depends only on the player's neighbors, the equilibrium (feedback) control depends on the entire network, requiring each player to look beyond its nearest neighbors or even its neighbors' neighbors. That said, we prove a correlation decay estimate (Proposition \ref{pr:corrdecay}), which shows that the covariance of two players' equilibrium state processes decays to zero with graph distance between these players.
Correlation decay is interesting in its own right, as it illustrates that asymptotic independence of players can arise both in a dense graph (because degrees are large) and a sparse graph (because typical vertices are very far apart).
In addition, we use correlation decay crucially in proving the convergence of the empirical distribution of state processes in equilibrium to a non-random limit (propagation of chaos), for large graph sequences $\{G_n\}$.

The key mathematical difficulty in the paper lies in the semi-explicit solution of the $n$-player game. As is standard for linear-quadratic $n$-player games, we reduce the problem to solving a coupled system of $n$ matrix differential equations of Riccati type. Riccati equations of this form do not often admit explicit solutions, but assuming the graph is transitive gives us enough symmetry to work with to derive a solution.

\section{Main results}

In this section we present all the main results of the paper, and we defer proofs to later sections. We first give the precise setup of the $n$-player game (Section \ref{se:finitehorizon}). After describing the semi-explicit solution of the equilibrium for transitive graphs (Section \ref{se:results_nplayer_finite}), we  then consider the large-$n$ behavior (Sections \ref{se:asymptotics} and \ref{se:examples}), paying particular attention to the distinction between the sparse and dense regimes.  Finally, we discuss the analogous cooperative game (Section \ref{se:cooperative}).

\subsection{The model setup.} \label{se:finitehorizon}

In this section we define a stochastic differential game associated to any finite graph $G = (V,E)$. All graphs will be simple and undirected. We abuse notation at times by identifying $G$ with its vertex set, e.g., by writing $v \in G$ instead of $v \in V$. Similarly, we write $|G|=|V|$ for the cardinality of the vertex set, and $\R^G=\R^V$ for the space of vectors indexed by the vertices.

Each vertex $v \in V$ is identified with a \emph{player}, and we associate to this player a state process on time horizon $T > 0$ with dynamics 
\begin{align}
dX^G_v(t) =  \alpha_v(t, \bm{X}^G(t)) dt + \sigma dW_v(t), \quad t \in [0,T],  \label{def:intro:SDE}
\end{align}
where $\sigma > 0$ is given, $(W_v)_{v \in V}$ are independent one-dimensional standard Brownian motions defined on a given filtered probability space $(\Omega, \F,\FF, \PP)$, and $\bm{X}^G = (X^G_v)_{v \in V}$ is the vector of state processes.  Players choose controls $\alpha_v$ from the set of (full-information) \emph{Markovian controls} $\A_G$, defined as the set of Borel-measurable functions $\alpha : [0,T] \times \R^V \to \R$ such that
\[
\sup_{(t,\bm x) \in [0,T] \times \R^V} \frac{|\alpha(t,\bm x)|}{1 + |\bm x|} < \infty.
\]
For any $\alpha_1,\ldots,\alpha_n \in \A_G$, the SDE system \eqref{def:intro:SDE} has a unique strong solution by a result of \cite{veretennikov1981strong} (see also \cite[Theorem 2.1]{krylov2005strong}).
The given initial states $\bm X^G(0)=(X^G_v(0))_{v \in V}$ are assumed non-random, and in many cases we will set them to zero for simplicity.

Each player $v \in V$ faces a quadratic cost function $J_v^G : \A_G^V \to \R$ that depends on the state processes of her nearest neighbors. For a non-isolated vertex $v$, we set
\begin{align}
J_v^G((\alpha_u)_{u \in V}) := \frac12\E\left[ \int_0^T |\alpha_v(t, \bm{X}^G(t))|^2dt + c\left|X^G_v(T)-\frac{1}{\deg_G(v)} \sum_{u \sim v} X^G_u(T) \right|^2 \right], \label{def:cost}
\end{align}
where $c > 0$ is a fixed constant, $\deg_G(v)$ denotes the degree of vertex $v$, and $u \sim v$ means that $(u,v)$ is an edge in $G$. For an isolated vertex $v$ (i.e., if $\deg_G(v)=0$), we set
\begin{align}
J_v^G((\alpha_u)_{u \in V}) := \frac12\E\left[ \int_0^T |\alpha_v(t, \bm{X}^G(t))|^2dt + c \left|X^G_v(T) \right|^2 \right]. \label{def:cost-isolated}
\end{align}

\begin{remark} \label{re:disconnected}
We gain little generality by allowing $G$ to be disconnected. Indeed, restricting attention to the connected components of $G$ yields decoupled games of the same form, which we can study separately. But when we discuss \Erdos \ and other random graphs, it is useful to fix a convention for how to handle isolated vertices. When discussing random graphs, we work always in the \emph{quenched} regime, with the realization of the graph frozen in the computation of the costs.
\end{remark}

In comparison to the usual settings of mean field games (MFGs), the key feature here is that the players do not interact with each other equally, but rather each player interacts (directly) only with her nearest neighbors in the graph. The form of the cost function implies indeed that each player, in addition to minimizing a standard quadratic energy term, will try to be as close as possible to the average of her nearest neighbors at the final time.
For this reason, we can think of this as a \emph{flocking model}.
The benchmark case to keep in mind is where $G$ is the complete graph on $n$ vertices, which corresponds to the usual MFG setup.

The first goal is to find a Markovian Nash equilibrium for this game, 
formally defined as follows, along with some generalizations. We write $\R_+=[0,\infty)$ throughout the paper.

\begin{definition} \label{def:Nash}
For a graph $G$ on vertex set $V=\{1,\ldots,n\}$ and a vector $\bm\epsilon=(\epsilon_i)_{i=1}^n \in \R_+^n$,  we say that a vector $\bm\alpha^* = (\alpha^{*}_i)_{i=1}^n \in \A_G^n$ of admissible strategies is a \emph{(Markovian) $\bm\epsilon$-Nash equilibrium on $G$} if 
\begin{equation*}
J_i^G(\bm\alpha^*) \le  \inf_{\alpha \in \A_G}J_i^G(\alpha^{*}_1, ..., \alpha^{*}_{i-1}, \alpha,  \alpha^{*}_{i+1}, ..., \alpha^{*}_n) + \epsilon_i, \quad \forall \, i=1,\ldots,n.
\end{equation*}
The corresponding \emph{equilibrium state process}  $\bm X^G=(X^G_i)_{i =1}^n$ is  the solution of the SDE
\begin{align*}
dX^G_i(t) = \alpha^*_i(t,\bm X^G(t))dt + \sigma dW_i(t).
\end{align*}
When the graph is understood from context, we may omit the qualifier ``on $G$."
When $\epsilon_1=\cdots=\epsilon_n=\epsilon$ for some $\epsilon \ge 0$, we refer to $\bm\alpha^*$ as a \emph{$\epsilon$-Nash equilibrium} instead of a \emph{$(\epsilon,\ldots,\epsilon)$-Nash equilibrium}.
Naturally, a \emph{$0$-Nash equilibrium} is simply called a  \emph{Nash equilibrium}. 
\end{definition}

The notion of $\epsilon$-Nash equilibrium for $\epsilon \ge 0$ is standard and means that no player can reduce her cost by more than $\epsilon$ by a unilateral change in control. The more general notion of $\bm\epsilon=(\epsilon_i)_{i=1}^n$-Nash equilibrium stated here is less standard, and it simply means that different players may stand to improve their costs by different amounts. Of course, a $\bm\epsilon=(\epsilon_i)_{i=1}^n$-Nash equilibrium is also a $\delta$-Nash equilibrium for $\delta=\max_{i=1}^n\epsilon_i$. But this distinction will be useful in asymptotic statements (as in the discussion after Theorem \ref{th:approxeq}), because the statement $\lim_{n\to\infty}\frac{1}{n}\sum_{i=1}^n\epsilon^n_i =0$ is of course much weaker than $\lim_{n\to\infty}\max_{i=1}^n\epsilon^n_i =0$ for triangular arrays $\{\epsilon^n_i : 1 \le i \le n\} \subset \R_+$.

\subsection{The equilibrium.} \label{se:results_nplayer_finite}

To solve the game described in Section \ref{se:finitehorizon}, we impose a symmetry assumption on the underlying graph.
Let $\Aut(G)$ denote the set of automorphisms of the graph $G=(V,E)$, i.e., bijections $\varphi : V \to V$ such that $(u,v) \in E$ if and only if $(\varphi(u),\varphi(v)) \in E$.  One should think of an automorphism as simply a \emph{relabeling} of the graph.

\begin{definition} \label{def:transitive}
We say $G$ is \emph{(vertex) transitive} if for every $u,v \in V$ there exists $\varphi \in \Aut(G)$ such that $\varphi(u)=v$.
\end{definition}

Essentially, a transitive graph ``looks the same" from the perspective of each vertex. Importantly, the game we are studying is clearly invariant under actions of $\Aut(G)$, in the sense that the equilibrium state process (if unique) should satisfy $(X^G_v)_{v \in V} \stackrel{d}{=} (X^G_{\varphi(v)})_{v \in V}$ for each $\varphi \in \Aut(G)$.
In the MFG setting, i.e., when $G$ is the complete graph, $\Aut(G)$ is the set of all permutations of the vertex set, and the $\Aut(G)$-invariance of the random vector $\bm X^G$  is better known as \emph{exchangeability}.
For a general graph, $\Aut(G)$ is merely a subgroup of the full permutation group, and we lose exchangeability.
While the transitivity of $G$ is a strong assumption, it is not surprising that a sufficiently rich group of symmetries would help us maintain some semblance of the tractability of MFG theory which stems from exchangeability.
Transitivity, in particular, ensures that we still have $X^G_v \stackrel{d}{=} X^G_u$ in equilibrium, for each $v,u \in V$.

We need the following notation.
Let $A_G$ denote the adjacency matrix of a graph $G$ on $n$ verties, and let $D_G = \mbox{diag}(\deg_G(1), ..., \deg_G(n))$ be the diagonal matrix of the degrees. If $G$ has no isolated vertices (i.e., all degrees are nonzero), we define the Laplacian\footnote{In the literature, there are several different matrices derived from a graph which go by the name \emph{Laplacian}. Our matrix $\lap_G$ is sometimes called the \emph{random walk normalized Laplacian} (or the negative thereof).} by
\begin{equation*}
\lap_G := D_G^{-1}A_G - I,
\end{equation*}
where $I$ is the identity matrix. It is easy to see that a transitive graph $G$ is always \emph{regular}, meaning each vertex has the same degree, which we denote $\delta(G)$. The Laplacian matrix then becomes $\lap_G = \tfrac{1}{\delta(G)}A_G-I$, which is notably a symmetric matrix.

\begin{remark} \label{re:spectrum}
Throughout the paper, we will make frequent use of the fact that $\lap_G$ has real eigenvalues, all of which are between $-2$ and $0$.
Indeed, note that $\lap_G = D_G^{-1/2}\widetilde{L}_G D_G^{1/2}$ where $\widetilde{L}_G=D_G^{-1/2}A_GD_G^{-1/2} - I$ is the \emph{symmetric normalized Laplacian}, and thus the eigenvalues of $\lap_G$ and $\widetilde{L}_G$ are the same; the properties of $\widetilde{L}_G$ are summarized by \cite[Sections 1.2 and 1.3]{chung-book}. Note that the all-ones vector is an eigenvector of $\lap_G$ with eigenvalue $0$.

\end{remark}

Our first main result is the following:

\begin{theorem}[Characterization of equilibrium on transitive graphs] \label{thm_equilibrium_n_player}
Suppose $G$ is a finite transitive graph on $n$ vertices without isolated vertices. Define $Q_G : \R_+ \to \R_+$ by
\begin{align}
Q_G(x) := (\det(I - x \lap_G))^{1/n}, \ \ \text{ for } x \in \R_+. \label{def:intro:Q}
\end{align}
Then $Q_G : \R_+ \to \R_+$ is well defined and continuously differentiable, and there exists a unique solution $f_G : [0,T] \to \R_+$ to the ODE
\begin{align*}
f_G'(t) = c Q_G'(f_G(t)), \hspace{0.5cm} f_G(0) = 0.
\end{align*}
Define $P_G : [0,T] \to \R^{n \times n}$ by
\begin{align}
P_G(t) := - f_G'(T-t) \lap_G \big(I - f_G(T-t) \lap_G \big)^{-1}, \label{def:intro:P*}
\end{align}
and finally define $\alpha^{G}_i \in \A_G$ for $i \in G$ by
\begin{equation*}
\alpha^{G}_i(t, \bm x) = - e_i^T P_G(t) \bm x,
\end{equation*}
where $(e_v)_{v \in G}$ is the standard Euclidean basis in $\R^G$. Then $(\alpha^{G}_i)_{i \in G}$ is a Nash equilibrium. For each $t \in (0,T]$, the equilibrium state process $\bm X^G(t)$ is normally distributed with mean vector $(I - f_G(T-t)\lap_G) (I  - f_G(T)\lap_G)^{-1}\bm X^G(0)$ and covariance matrix
\begin{align*}
\sigma^2 (I - f_G(T-t)\lap_G)^2 \int_0^t (I - f_G(T-s)\lap_G)^{-2}ds.
\end{align*}
Finally, writing $|\cdot|$ for the Euclidean norm, the time-zero average value is
\begin{align}
\mathrm{Val}(G) := \frac{1}{n}\sum_{v \in G} J^G_v((\alpha^{G}_i)_{i \in G}) = \frac{|P_G(0)\bm X^G(0)|^2}{2\tr(P_G(0))} - \frac{\sigma^2}{2} \log \frac{\tr(P_G(0))}{nf'_G(T)}.		\label{def:EQvalue}
\end{align}
\end{theorem}

The proof is given in Section \ref{section:n_player_finite}.
As usual, we first reduce our (linear-quadratic) game to a system of matrix differential equations of Riccati type in Section \ref{se:HJBs}. In our setting we can explicitly solve these Riccati equations using symmetry arguments based on the transitivity assumption.
In Section \ref{se:extension_general_matrices}, we discuss an extension of Theorem \ref{thm_equilibrium_n_player} to a more general class of matrices $L$, or equivalently to weighted graphs, satisfying a suitable generalization of the transitivity assumption.

It is important to note that the equilibrium controls $\alpha^G_i$ obtained in Theorem \ref{thm_equilibrium_n_player} are \emph{nonlocal}, in the sense that the control of player $i$ depends on the states of all of the players, not just the neighbors. Naive intuition would suggest that players should only look at the states of their neighbors, because the objective of each player is to align at time $T$ with those neighbors. On the contrary, a rational player anticipates that her neighbors will in turn try to align with their own neighbors, which leads the player to follow the states of the neighbors' neighbors, and similarly the neighbors' neighbors' neighbors, and so on.

It is worth noting that in the setting of Theorem \ref{thm_equilibrium_n_player} we have
\begin{align*}
\E\left[\frac{1}{n}\sum_{v \in G}X^G_v(t)\right] = \frac{1}{n}\sum_{v \in G}X^G_v(0).
\end{align*}
That is, the average location of the players stays constant over time, in equilibrium.
Indeed, this follows easily from the formula for the mean $\E[\bm X^G(t)]$ and from the fact that the vector of all ones is an eigenvector with eigenvalue 0 for the symmetric matrix $\lap_G$.

We suspect that the Markovian Nash equilibrium identified in Theorem \ref{thm_equilibrium_n_player} is the unique one. This could likely be proven using similar arguments to those of \cite[Section II.6.3.1]{carmona-delarue-book}, but for the sake of brevity we do not attempt to do so.

\subsection{Asymptotic regimes.} \label{se:asymptotics}

The form of the equilibrium computed in Theorem \ref{thm_equilibrium_n_player} lends itself well to large-$n$ asymptotics after a couple of observations. First, for simplicity, we focus on the case $\bm X^G(0)=\bm{0}$. Transitivity of the graph $G$ (or Lemma \ref{le:symmetries}) ensures that $X^G_i(t) \stackrel{d}{=} X^G_j(t)$ for all $i,j \in G$ and $t > 0$, and we deduce that each $X^G_i(t)$ is a centered Gaussian with variance
\begin{align}
\Var( X^G_i(t)) &= \frac{1}{n}\sum_{k=1}^n \Var( X^G_k(t) ) \nonumber \\
	&= \frac{\sigma^2}{n}\tr\left[ (I-f_G(T-t)\lap_G)^2 \int_0^t (I-f_G(T-s)\lap_G)^{-2}ds\right] \nonumber \\
	&=  \frac{\sigma^2}{n} \sum_{k=1}^n \int_0^t  \left(\frac{1 - f_G(T-t)\lambda^G_k}{1 - f_G(T-s)\lambda^G_k}\right)^2 ds, \label{intro:variance-lambda}
\end{align}
where $\lambda_1^G,\ldots,\lambda_n^G$ are the eigenvalues of $\lap_G$, repeated by multiplicity.
The average over $k=1,\ldots,n$ can be written as an integral with respect to the empirical eigenvalue distribution,
\begin{align}
\mu_G := \frac{1}{n}\sum_{i=1}^n\delta_{\lambda^G_i}, \label{def:spectralmeasure}
\end{align}
which we recall is supported on $[-2,0]$, as in Remark \ref{re:spectrum}.
The other quantities in Theorem \ref{thm_equilibrium_n_player} can also be expressed in terms of $\mu_G$. Indeed, the value $\mathrm{Val}(G)$ becomes
\begin{equation}
\mbox{Val}(G) = - \frac{\sigma^2}{2} \log \frac{\tr{(P_{G}(0))}}{n f'_{G}(T)} = - \frac{\sigma^2}{2} \log  \int_{[-2,0]} \frac{- \lambda}{1 - f_{G}(T)\lambda} \mu_{G}(d\lambda), \label{def:valueidentity}
\end{equation}
and the function $Q_G$ defined in \eqref{def:intro:Q} becomes
\begin{align*}
Q_G(x) &= \left(\prod_{i=1}^n(1-x\lambda^G_i)\right)^{1/n} = \exp \int_{[-2,0]} \log(1-x\lambda) \,\mu_G(d\lambda).
\end{align*}
Thus, if we are given a sequence of graphs $G_n$ such that $\mu_{G_n}$ converges weakly to some probability measure, it is natural to expect the equilibrium computed in Theorem \ref{thm_equilibrium_n_player} to converge in some sense. This is the content of our next main result, which we prove in Section \ref{se:asymptoticregimes}:

\begin{theorem}[Large-scale asymptotics on transitive graphs] \label{thm:convergence_X^Gn_k}
Let $\{G_n\}$ be a sequence of finite transitive graphs without isolated vertices, with $\lim_{n\to\infty}|G_n|=\infty$. Let $\bm X^{G_n}$ denote the equilibrium state process identified in Theorem \ref{thm_equilibrium_n_player}, started from initial position $\bm X^{G_n}(0)= \bm{0}$.
Suppose $\mu_{G_n}$ converges weakly to a probability measure $\mu$, and define $Q_\mu : \R_+ \to \R_+$ by
\begin{align*}
Q_\mu(x) := \exp \int_{[-2,0]} \log(1-x\lambda) \,\mu(d\lambda).
\end{align*}
Then the following holds:
\begin{enumerate}
\item There exists a unique solution $f_\mu : [0,T] \to \R_+$ of the ODE
\begin{align}
f'_\mu(t) = cQ'_\mu(f_\mu(t)), \quad f_\mu(0)=0. \label{def:limitODE}
\end{align}
\item For any vertex sequence $k_n \in G_n$ and any $t \in [0,T]$, the law of $X^{G_n}_{k_n}(t)$ converges weakly as $n\to\infty$ to the Gaussian distribution with mean zero and variance
\begin{align}
\varfunc_\mu(t) =  \sigma^2\int_0^t\int_{[-2,0]} \left(\frac{1 - \lambda f_\mu(T-t)}{1 - \lambda f_\mu(T-s)}\right)^2\mu(d\lambda) ds. \label{def:gaussianlimit}
\end{align}
\item For any $t \in [0,T]$, the (random) empirical measure $\frac{1}{|G_n|}\sum_{i \in G_n} \delta_{X^{G_n}_i(t)}$ converges weakly in probability as $n\to\infty$ to the (non-random) Gaussian distribution $\mathcal{N}(0,\varfunc_\mu(t))$.
\item The time-zero values given in \eqref{def:EQvalue} with $\bm X^{G_n}(0)= 0$ converge:
\begin{align}
\lim_{n\to\infty}\mathrm{Val}(G_n) = - \frac{\sigma^2}{2} \log \int_{[-2,0]} \frac{-\lambda}{1-\lambda f_\mu(T)}\mu(d\lambda). \label{def:limitvalue}
\end{align}
\end{enumerate}
\end{theorem}

There are many concrete graph sequences $\{G_n\}$ for which $\mu_{G_n}$ can be shown to converge to a tractable (typically continuous) limiting measure, and we document several notable cases in Section \ref{se:examples}. The Laplacian spectrum is in fact quite tractable and well-studied.
There is a substantial literature on the eigenvalues of Laplacian (and other) matrices of graphs (\cite{chung-book,godsil2013algebraic}), which are well known to encode significant structural information about the graph. 

In addition, graph convergence concepts like \emph{local weak convergence} are known to imply weak convergence of the spectral measure (\cite{BordenaveLelarge}); see Section \ref{se:intrinsic} for some further discussion.

\begin{remark} \label{re:variance}
We develop in Section \ref{se:ODEanalysis} some noteworthy qualitative and quantitative properties of the equilibrium variance $\varfunc_\mu(t)$ given in \eqref{def:gaussianlimit}.
We show in Proposition \ref{pr:fQproperties} that  $\varfunc_\mu(0)=0$, $\varfunc'_\mu(0)=\sigma^2$, and $\varfunc''_\mu(0)=-2\sigma^2 (f'_\mu(T))^2/cQ_\mu(f_\mu(T))$. In particular, for short times, the leading-order behavior $\varfunc_\mu(t) = \sigma^2 t + o(t)$ does not depend on the underlying graph. It is only at the second order or at longer time horizons that the influence of the graph is felt.
\end{remark}

\begin{remark}
The restriction to $\bm X^{G_n}(0)=\bm 0$ in Theorem \ref{thm:convergence_X^Gn_k} is merely to simplify the resulting formulas. One could easily accommodate the more general setting in which the empirical measure of initial states converges to some limiting distribution. 
(Note, however, that if the graph is not transitive, then general initial states may confound the convergence analysis; see  \cite{coppini2019law}, and especially Remark 1.2 therein for a relevant discussion of uncontrolled models.)
 In addition, a functional version of Theorem \ref{thm:convergence_X^Gn_k} can likely be derived under no further assumptions, in which the Gaussian process $(X^{G_n}_{k_n}(t))_{t \in [0,T]}$ converges weakly in $C([0,T])$ to a limiting Gaussian process. We omit these generalizations, as the more complicated statements do not shed any light on the role of the network structure, which is the main focus of our work.
\end{remark}

\subsubsection{Dense graphs.}
If $G_n$ is the complete graph, then it turns out that $\mu_{G_n} = \frac{1}{n}\delta_0 + \frac{n-1}{n}\delta_{-n/(n-1)} \to \delta_{-1}$, which leads to a simpler form for the limiting law in \eqref{def:gaussianlimit}. More generally, the case $\mu_{G_n} \to \delta_{-1}$ represents a ``dense" regime, as described in the following result. Recall that all transitive graphs are \emph{regular}, meaning each vertex has the same degree. The following is proven in Section \ref{se:examples-proofs}:

\begin{corollary}[Large-scale asymptotics on dense transitive graphs] \label{co:dense}
Suppose $\{G_n\}$ is a sequence of transitive graphs, and suppose each vertex of $G_n$ has common degree $\delta(G_n) \ge 1$.
Then $\mu_{G_n} \rightarrow \delta_{-1}$ if and only if $\delta(G_n)\to \infty$. In this case, the limiting variance \eqref{def:gaussianlimit} and value \eqref{def:limitvalue} simplify to
\begin{align}
\varfunc_{\delta_{-1}}(t) = \sigma^2 t \frac{ 1 + c(T-t) }{1+cT}, \qquad
\lim_{n\to\infty}\mathrm{Val}(G_n) =  \frac{\sigma^2}{2} \log(1+ cT).  \label{def:vardense}
\end{align}
Moreover, there is a constant $C < \infty$, depending only on $c$ and $T$, such that
\begin{align}
|\varfunc_{G_n}(t) - \varfunc_{\delta_{-1}}(t)| + \left|\mathrm{Val}(G_n) - \tfrac{\sigma^2}{2} \log(1+ cT)\right| \le C/\delta(G_n), \quad 
\forall n \in \N, \ t \in [0,T]. \label{def:ineq}
\end{align}
Finally, the Gaussian law $\mathcal{N}(0,\varfunc_{\delta_{-1}}(t))$ is precisely the time-$t$ law of the unique solution  of the SDE
\begin{align}
dX(t) = -\frac{cX(t)}{1+c(T-t)}dt + \sigma dW(t), \quad X(0)=0. \label{def:denseSDE}
\end{align}\\
\end{corollary}

Corollary \ref{co:dense} shows that the dense regime is particularly tractable. In particular, the mean field case (where $G_n$ is the complete graph) is \emph{universal} in the sense that the same limit arises for any other transitive graph sequence with diverging degree. Moreover, the rate $C/\delta(G_n)$ in \eqref{def:ineq} becomes $C/n$ in the mean field case, which is the best-known convergence rate for the value functions of well-behaved MFGs (\cite[Theorem 2.13]{CDLL}).

\begin{remark} \label{re:variance-dense}
We show in Proposition \ref{pr:vardeltamin} that the dense graph regime uniquely achieves the lowest possible variance; precisely, we have $\varfunc_\mu(t) \ge \varfunc_{\delta_{-1}}(t)$, recalling the definitions  \eqref{def:gaussianlimit} and \eqref{def:vardense}, with equality only when $\mu=\delta_{-1}$.
The example of the torus graphs in Section \ref{se:ex:torus} below illustrates what appears to be a general principle, that \emph{a highly connected graph has smaller variance in equilibrium}. This makes intuitive sense, as a higher degree means each player has a larger set of neighbors to be attracted toward. 
\end{remark}

Our next result shows that in the dense regime we may use the limiting object to construct approximate equilibria for $n$-player games on general large dense graphs (not necessarily transitive), in the same way that the equilibrium of a MFG can be used to build approximate equilibria for finite games.

\begin{theorem}[Approximate equilibria on general dense graphs] \label{th:approxeq}
Suppose $G$ is a finite graph. For each vertex $v$ of $G$, define a control
\begin{align*}
\alpha_v^{\mathrm{MF}}(t,\bm x) := \frac{- cx_v}{1+c(T-t)}, \quad t \in [0,T], \ \ \bm x = (x_u)_{u \in G} \in \R^G.
\end{align*}
Finally, define $\bm\epsilon^G=(\epsilon^G_v)_{v \in G} \in \R_+^G$ by 
\begin{align*}
\epsilon^G_v := \left\{ \begin{array}{ll}
    \sigma^2\frac{cT}{1+cT}\sqrt{\frac{cT(2+cT)}{ \deg_G(v)}} &  \mbox{ if } \deg_G(v) \ge 1 \\
    0 & \mbox{ if } \deg_G(v) = 0.
\end{array}
\right.
\end{align*}
Then, for each $n$, $(\alpha_v^{\mathrm{MF}})_{v \in G}$ is an $\bm\epsilon^G$-Nash equilibrium on $G$.
In particular, if\footnote{As usual, we write $a \vee b := \max\{a,b\}$.}
\begin{align*}
\epsilon_G := \sigma^2\frac{cT}{1+cT}\sqrt{\frac{cT(2+cT)}{1 \vee \delta(G)}}, \qquad \text{where } \ \ \delta(G):= \min_{v \in G}\deg_G(v),
\end{align*}
then $(\alpha_v^{\mathrm{MF}})_{v \in G}$ is a $\epsilon_G$-Nash equilibrium on $G$.
\end{theorem}

We use the notation $\alpha_v^{\mathrm{MF}}$ because this is precisely the control one obtains from the corresponding MFG (see Lemma \ref{lemma_alpha*_optimal}).

Hence, Theorem \ref{th:approxeq} says that on a graph sequence with ``diverging degree" in some sense, the MFG provides a (decentralized, symmetric) approximate Nash equilibrium. 
More precisely, if $\{G_n\}$ is a sequence of graphs with diverging minimal degree $\delta(G_n) \to \infty$, then the controls $\bm\alpha^n:=(\alpha_v^{\mathrm{MF}})_{v \in G_n}$ form an $\epsilon_{G_n}$-Nash equilibrium for each $n$ with $\lim_n\epsilon_{G_n} = 0$.
Of course, the now-classical theory of MFGs tells us the same thing when $G_n$ is the complete graph (see, e.g., \cite[Section II.6.1]{carmona-delarue-book}, or \cite[Theorem 12]{huang2006large} for the standard rate of $\epsilon_{G_n} = O(1/\sqrt{n})$), but Theorem \ref{th:approxeq} gives a threshold of \emph{how dense} the graph needs to be in order for the mean field approximation to remain valid.
The constant $\epsilon_{G_n}$ shows quantitatively how the accuracy of the mean field approximation depends on the ``denseness" of the graph, as measured by the minimal degree.
Some examples beyond the complete graph will be discussed in Section \ref{se:examples} below.

In fact, we may relax the denseness threshold if we are happy to assert that $(\alpha_v^{\mathrm{MF}})_{v \in G_n}$ form an approximate equilibrium in a weaker sense, suggested by the most general form of Definition \ref{def:Nash}. A small fraction of players (namely, those with small degree) might have a lot to gain by deviating, but this potential gain from deviation is small when averaged over all players. Precisely, suppose that instead of the minimum degree diverging, we suppose merely that degrees diverge in the following averaged sense:
\begin{align}
\lim_{n\to\infty}\frac{1}{|G_n|}\sum_{v \in G_n}(1 \vee \deg_{G_n}(v))^{-1/2} = 0. \label{def:avgdegreediverge}
\end{align}
Then $(\alpha_v^{\mathrm{MF}})_{v \in G_n}$ is an $\bm\epsilon^{G_n}$-Nash equilibrium, and $\lim_n\frac{1}{n}\sum_{i=1}^n\epsilon^{G_n}_i = 0$.

In summary, different manners of quantifying the concept of \emph{approximate equilibrium} lead to different sparsity/denseness thresholds for the validity of the mean field approximation. The \Erdos \ case in Section \ref{se:ERexample} gives a concrete example.

\subsubsection{Correlation decay and asymptotic independence.}
Before discussing examples, we lastly present an estimate of correlation decay, which is crucial in the proof of convergence of the empirical measure in Theorem \ref{thm:convergence_X^Gn_k}, and which also reveals what form of \emph{asymptotic independence} between the players can be expected. See Section \ref{subsection:cov_bounds} for the proof:

\begin{proposition}[Correlation decay on transitive graphs] \label{pr:corrdecay}
Let $G$ be a finite transitive graph without isolated vertices, and let $\bm X^G$ denote the equilibrium state process identified in Theorem \ref{thm_equilibrium_n_player}.
Suppose each vertex of $G$ has degree $\delta(G) \in \N$.
For vertices $u,v \in G$, let $d_G(u,v)$ denote the graph distance, defined as the length of the shortest path from $u$ to $v$ (and $\infty$ if no such path exists). Let $\gamma = cT/(1+cT) \in (0,1)$. Then
\begin{align}
|\Cov(X^G_u(t), X^G_v(t))| \le 2\sigma^2 t\frac{\gamma^{d_G(u,v)} \big(1 + d_G(u,v)(1 - \gamma)\big)}{\delta(G)(1 - \gamma)^2}1_{\{d_G(u,v) < \infty\}}. \label{def:corrdecaybound}
\end{align}
\end{proposition}

Note that the right-hand side of \eqref{def:corrdecaybound} is a bounded function of $\delta(G)$ and $d_G(u,v)$.
If $G_n$ is the complete graph on $n$ vertices (i.e., the mean field case), then $\delta(G_n)=n-1\to\infty$, and each pair of players (in equilibrium) becomes asymptotically independent as $n\to\infty$. This is an instance of the phenomenon of \emph{propagation of chaos} for mean field systems. More generally, this remains true for any dense graph sequence, i.e., whenever $\delta(G_n) \to \infty$.

On the other hand, the picture is rather different for a sparse graphs sequence, i.e., when $\sup_n\delta(G_n) < \infty$. An \emph{arbitrary} pair of players can no longer be expected to become asymptotically independent as $n\to\infty$, but only \emph{distant} players. More precisely, two players $u_n,v_n \in G_n$ become asymptotically independent only if $d_{G_n}(u_n,v_n) \to \infty$ (since $\gamma < 1$). For a transitive graph sequence with $|G_n| \to \infty$ and $\sup_n\delta(G_n)  < \infty$, it is always the case that the distance between two \emph{uniformly random} vertices converges to infinity in probability, and it follows that two (uniformly) randomly chosen players are asymptotically independent.

In summary, for a sequence of (transitive) graphs $G_n$ with $|G_n|\to\infty$, asymptotic independence of a typical pair of players arises for one of two quite distinct reasons. Either:
\begin{enumerate}
\item The degree diverges, and each player interacts with many other players, with pairwise interaction strengths of order $1/\delta(G_n) \to 0$.
\item The degrees stay bounded, but typical players are very far apart in the graph and thus very weakly correlated.
\end{enumerate}

The correlation decay estimate of Proposition \ref{pr:corrdecay} is the key ingredient which allows us to deduce Theorem \ref{thm:convergence_X^Gn_k}(iii) from Theorem \ref{thm:convergence_X^Gn_k}(ii), i.e., to prove the empirical measure convergence. Indeed, we will use this covariance bound along with the Gaussian Poincar\'e inequality to prove that $X^{G_n}_{v_n}(t)$ and $X^{G_n}_{u_n}(t)$ are asymptotically independent as $n\to\infty$, when $u_n$ and $v_n$ are independent uniformly random vertices in $G_n$ for each $n$. That is, $(X^{G_n}_{v_n}(t),X^{G_n}_{u_n}(t))$ converges in joint law to $\mathcal{N}(0,\varfunc_\mu(t))^{\otimes 2}$.  By a standard propagation of chaos argument, this is equivalent to the convergence of the empirical measure $\frac{1}{|G_n|}\sum_{i \in G_n} \delta_{X^{G_n}_i(t)}$ to $\mathcal{N}(0,\varfunc_\mu(t))$. See Section \ref{se:convempmeas} for details. 

Moreover, the empirical measure convergence is also equivalent to the convergence in joint law of $(X^{G_n}_{v_n^i}(t))_{i=1}^k$ to  $\mathcal{N}(0,\varfunc_\mu(t))^{\otimes k}$ as $n\to\infty$ for fixed $k \ge 2$, where $(v^1_n,\ldots,v^k_n)$ are either independent uniformly random vertices or a uniformly random choice from the $n(n-1)\cdots(n-k+1)$ possible $k$-tuples of distinct vertices. Thus, in both the sparse and dense regime, we obtain a full picture of propagation of chaos, up to a randomization of the choice of vertices.

\subsection{Examples.} \label{se:examples}

In this section we specialize the results of Section \ref{se:asymptotics} to a short (and by no means exhaustive) list of somewhat tractable natural large-graph models. We focus on cases where the minimum degree and/or the empirical eigenvalue distribution of the graph are tractable, as these quantities are particularly relevant to the main results of Section \ref{se:asymptotics}.

\subsubsection{The complete graph.}
Let us summarize what we have mentioned regarding the simplest (mean field) case, where $G_n$ is the complete graph on $n$ vertices. In this case, the Laplacian matrix takes the form
\begin{align*}
\lap_{G_n}=\frac{1}{n-1}(J-I) - I,
\end{align*}
where $J$ is the matrix of all ones. From this we easily deduce that the eigenvalues of $\lap_{G_n}$ are $0$ and $-\tfrac{n}{n-1}$, with respective multiplicities $1$ and $n-1$. Hence, $\mu_{G_n} \to \delta_{-1}$, and the degree $\delta(G_n)=n-1 \to \infty$. The complete graph is of course transitive, and all of our main theorems apply, in particular Corollary \ref{co:dense}.

In the complete graph setting, our model essentially becomes the $\epsilon=q=0$ case of \cite{carmona2013mean}. The only difference is that in \cite{carmona2013mean} each player is included in the empirical average; that is, the terminal cost of player $k$ is $|\frac{1}{n}\sum x_i - x_k|^2$ instead of $|\frac{1}{n-1}\sum_{i\neq k}x_i - x_k|^2$.
This can easily be fit into our framework, as the following remark explains.

\begin{remark} \label{re:self-loops}
For a finite transitive graph $G$ without isolated vertices, every vertex has a common degree $\delta(G)$. Letting $N_v(G)$ denote the union of $\{v\}$ and the set of neighbors of a vertex $v$ in $G$, we can write the terminal cost function for player $v$ as
\begin{align*}
\left|\frac{1}{\delta(G)}\sum_{u \sim v}x_u - x_v \right|^2 = \left(\frac{\delta(G)+1}{\delta(G)} \right)^2\left|\frac{1}{\delta(G)+1}\sum_{u \in N_v(G)}x_u - x_v \right|^2.
\end{align*}
Hence, we can modify our setup so that each each player is included in the average in the terminal cost, simply by modifying the constant $c$ by a factor of $(1+1/\delta(G))^2$.
\end{remark}

\subsubsection{The cycle graph.}

Suppose now that $G_n=C_n$ is the cycle on $n$ vertices. This is a transitive graph in which every vertex has common degree $\delta(C_n)=2$. The adjacency matrix $A_{C_n}$ is a circulant matrix, which makes it easy to calculate the eigenvalues as $2 \cos( 2\pi k/n)$ for $k=1,\ldots,n$. The eigenvalues of the Laplacian $\lap_{C_n} = \tfrac12 A_{C_n}-I$ are thus $\lambda^{C_n}_k=\cos(2\pi k/n) - 1$ for $k=1,\ldots,n$. In this case, for a bounded continuous function $f$ we compute
\begin{align*}
\int f\,d\mu_{C_n} = \frac{1}{n}\sum_{k=1}^n f(\cos(2\pi k/n) - 1) \to \int_0^1 f(\cos(2\pi u) - 1) \,du, \ \text{ as } n \to \infty,
\end{align*}
which shows that $\mu_{C_n}$ converges weakly to the probability measure $\mu$ given by the law of $\cos(2\pi U) - 1$, where $U$ is uniform in $[0,1]$, i.e., $\mu(dx)=1_{[-2,0]}(x)\frac{dx}{\pi\sqrt{-x(2+ x)}}$. The function $Q_\mu$ in Theorem \ref{thm:convergence_X^Gn_k} is then
\begin{align}
Q_\mu(x) = \exp \int_0^1 \log\big(1 + x - x\cos(2\pi u)\big) du, \qquad x \ge 0. \label{def:Qcycle}
\end{align}
In Section \ref{se:cyclegraphproof} we derive a semi-explicit solution of the ODE \eqref{def:limitODE} in this setting:\\

\begin{proposition} \label{prop:cyclegraph}
Define $Q_\mu$ as in \eqref{def:Qcycle}. Then $Q_\mu(x)=\tfrac12(\sqrt{1+2x}+x+1)$, and the unique solution of the ODE $f'_\mu(t)=cQ'_\mu(f_\mu(t))$ with $f_\mu(0)=0$ is given by
\begin{align*}
f_\mu(t) = \Phi^{-1}\left(\log 2 +\frac{ct-1}{2} \right),
\end{align*}
where $\Phi^{-1}$ is the inverse of the strictly increasing function $\Phi : \R_+\to\R_+$ defined by
\begin{align*}
\Phi(x) := \log(1+\sqrt{1+2x}) - \sqrt{1+2x} + x + \tfrac12.
\end{align*}
\end{proposition}

\begin{figure}[h]
    \centering
    \includegraphics[width=0.8\columnwidth]{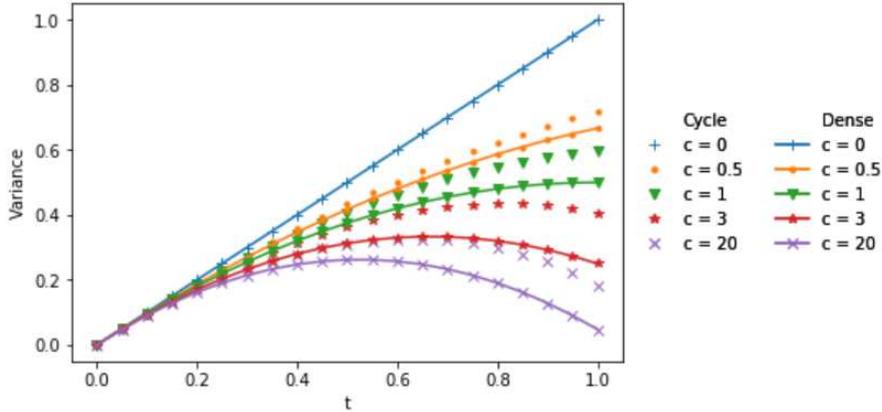}
    \caption{Variance of a typical player over time in the dense graph (line with markers, equation \eqref{def:vardense}) and in the cycle graph (markers only, equation \eqref{def:cyclevariance}) for different values of $c$. Here $T=\sigma=1$.}
    \label{fig:cycleVSdense}
\end{figure}

The variance from \eqref{def:gaussianlimit} then becomes
\begin{align}
\varfunc_\mu(t) =  \sigma^2\int_0^t\int_0^1 \left(\frac{1 - (\cos(2\pi u) -1) \Phi^{-1}\left(\log 2 + \frac{c(T-t)-1}{2} \right)}{1 - (\cos(2\pi u) -1) \Phi^{-1}\left(\log 2 + \frac{c(T-s)-1}{2} \right)}\right)^2 du \, ds. \label{def:cyclevariance}
\end{align}
This does not appear to simplify further, but Figure \ref{fig:cycleVSdense} gives plots for various $c$.
Note that the variance in the dense case is always lower than in the cycle case, as we show in Proposition \ref{pr:fQproperties}. In both cases, the variance at any fixed time $t$ decreases with $c$. 

As $c\to\infty$, the variance $\varfunc_\mu(t)$ in both the dense and cycle graph cases can be shown to converge to $\sigma^2 t(T-t)/T$, which is the same as that of a Brownian bridge.

\subsubsection{The torus.} \label{se:ex:torus}

For $d \in \N$, consider the torus $G_n=\Z^d_n:= \Z^d/ n \Z^d$. That is, this graph is the subgraph of the integer lattice $\Z^d$ with vertex set $\{1,\ldots,n^d\}$ and with ``wrapping around" at the boundary. 
The eigenvalues of $\lap_{\Z^d_n}$ are easily computed from those of $\lap_{C_n}$, the cycle graph from the previous section, after noting that $\Z^d_n$ is the $d$-fold Cartesian product of the cycle $C_n$ with itself.
In particular, if $G$ and $H$ are two graphs, and $A_G$ and $A_H$ have eigenvalues $(\eta^G_v)_{v \in G}$ and $(\eta^H_v)_{v \in H}$ respectively, then the eigenvalues of the adjacency matrix of the Cartesian product of $G$ and $H$ are given by $(\eta^G_v+\eta^H_u)_{u \in G, v \in H}$.\footnote{See Chapter 7.14 of \cite{godsil2013algebraic} for definition of the Cartesian product of graphs and Chapter 9.7 for a derivation of the eigenvalues of the adjacency matrix of a Cartesian product.} In particular, the eigenvalues of $A_{\Z^d_n}$ are
\begin{align*}
\sum_{i=1}^d 2\cos(2\pi k_i/n), \qquad \bm k=(k_1,\ldots,k_d) \in \Z^d_n.
\end{align*}
Noting that each vertex in $\Z^d_n$ has degree $2d$, we find that the eigenvalues of $\lap_{\Z^d_n}=\tfrac{1}{2d}A_{\Z^d_n} - I$ are
\begin{align*}
\lambda^{\Z^d_n}_{\bm k} = \frac{1}{d}\sum_{i=1}^d \cos(2\pi k_i/n) - 1, \qquad \bm k=(k_1,\ldots,k_d) \in \Z^d_n.
\end{align*}
Hence, for a bounded continuous function $f$ we compute
\begin{align*}
\int f\,d\mu_{\Z^d_n} &= \frac{1}{|\Z^d_n|}\sum_{\bm k \in \Z^d_n} f_\mu(\lambda^{\Z^d_n}_{\bm k}) = \frac{1}{n^d}\sum_{k_1,\ldots,k_d=1}^n f\left(\frac{1}{d}\sum_{i=1}^d \cos(2\pi k_i/n) - 1\right) \\
	&\to \int_{[0,1]^d}f\left(\frac{1}{d}\sum_{i=1}^d \cos(2\pi u_i) - 1\right)du, \qquad \text{ as } \ \ n\to\infty,
\end{align*}
which shows that $\mu_{\Z^d_n}$ converges weakly to the probability measure $\mu$ given by the law of $\frac{1}{d}\sum_{i=1}^d \cos(2\pi U_i) - 1$, where $U_1,\ldots,U_d$ are independent uniform random variables in $[0,1]$. The function $Q_\mu$ in Theorem \ref{thm:convergence_X^Gn_k} is then
\begin{align}
Q_\mu(x) = \exp \int_{[0,1]^d} \log\left(1 + x - \frac{x}{d}\sum_{i=1}^d \cos(2\pi u_i)\right) du. \label{def:Qtorus}
\end{align}

\begin{figure}[h]
    \centering
    \includegraphics[width=0.6\columnwidth]{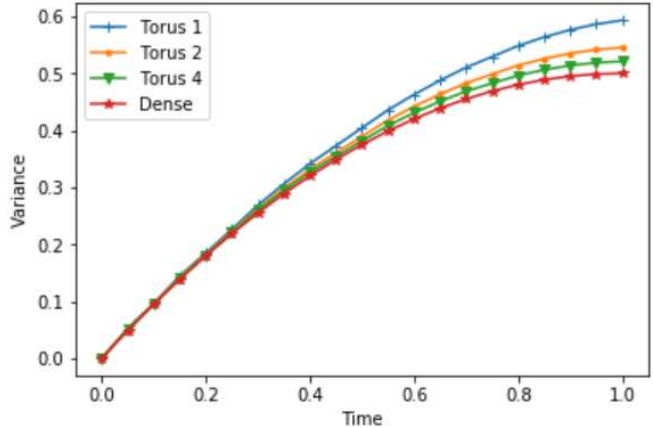}
    \caption{Variance of a typical player over time in the torus of dimensions $d=1,2,4$ and the dense case. Here $T=c=\sigma=1$.}
    \label{fig:torus}
\end{figure}

We cannot evaluate \eqref{def:Qtorus} or the solution $f$ of the ODE \eqref{def:limitODE} explicitly, for the torus of dimension $d > 1$ (the case $d=1$ is the cycle graph). But we can easily do so numerically. Figure \ref{fig:torus} shows the variance $\varfunc_\mu(t)$ of \eqref{def:gaussianlimit} for the torus of various dimensions, compared with the dense graph case. Notably, the variance decreases with the dimension $d$, supporting the intuition that a more highly connected graph leads to a behavior closer to the mean field regime.

\subsubsection{\Erdos \ graphs.} \label{se:ERexample}

Most of our main results require a transitive graph and thus have little to say about classical random graph models, such as \Erdos, random regular graphs, or the configuration model, which generate graphs which are non-transitive with high probability. In particular, in applications of Theorems \ref{thm_equilibrium_n_player} and \ref{thm:convergence_X^Gn_k} we cannot take any advantage of the vast body of literature on the behavior of the eigenvalue distribution of the adjacency and Laplacian matrices of these (non-transitive) random graph models.
That said, we mention here some noteworthy \emph{dense} random graph models, to which Theorem \ref{th:approxeq} applies and shows that the MFG approximation is valid.

For the \Erdos \ graph $G(n,p_n)$, it is known that as $n\to\infty$ with $\liminf_n np_n/\log n > 1$,  the minimal degree converges to infinity in probability (\cite[Lemma 6.5.2]{durrett2007random}). Hence, Theorem \ref{th:approxeq} applies in this regime to give a random sequence $\epsilon_n \ge 0$ converging to zero in probability such that $(\alpha^{\mathrm{MF}}_v)_{v \in G_n}$ is an $\epsilon_n$-Nash equilibrium for each $n$. This is sharp in a sense, because $p_n >  \log n /n$ is precisely the threshold for connectedness: If $\limsup_n np_n/\log n < 1$, then $G(n,p_n)$ contains isolated vertices with high probability (in particular, the minimal degree is $1$), and we cannot expect $\alpha^{\mathrm{MF}}_v$ to be near-optimal for $v$ in a small connected component. This might be compared to the main result of \cite{delarue2017mean}, which  keeps $p_n=p$ constant as $n\to\infty$ and finds likewise that the usual MFG approximation is valid for a class of games on the \Erdos \ graph.

If we relax our concept of \emph{approximate equilibrium}, as in the discussion after Theorem \ref{th:approxeq}, then we may push the denseness threshold all the way to $np_n\to\infty$. That is, if $np_n\to\infty$, then a straightforward calculation using the fact that $\deg_{G_n}(v) \sim \mathrm{Binomial}(n-1,p_n)$ shows that 
\begin{align*}
\E\left[\frac{1}{n}\sum_{v =1}^n (1 \vee \deg_{G_n}(v))^{-1/2}\right] = \sum_{k=0}^{n-1}{\binom{n-1}{k}}p_n^k(1-p_n)^{n-k-1}(1 \vee k)^{-1/2} \to 0,
\end{align*}
which ensures by Theorem \ref{th:approxeq} that there exist random (graph-dependent) variables $\bm\epsilon^n=(\epsilon^n_v)_{v=1}^n$ such that $\frac{1}{n}\sum_{v=1}^n\epsilon^n_v \to 0$ in probability and $(\alpha^{\mathrm{MF}}_v)_{v \in G_n}$ forms a $\bm\epsilon^n$-Nash equilibrium.
Note that this threshold $np_n\to\infty$ means that the expected degree (of a randomly chosen vertex) diverges.

In the extremely sparse (diluted) regime $\lim_n np_n = \theta \in (0,\infty)$, the degree of a typical vertex converges in law to Poisson($\theta$), and Theorem \ref{th:approxeq} yields no information.

These thresholds are in line with recent work on interacting particle systems (without game or control). For example, interacting particle systems on either the complete graph or the \Erdos \ graph $G(n,p_n)$ converge to the same mean field (McKean-Vlasov) limit as soon as $np_n\to\infty$  (\cite{bhamidi2019weakly,oliveira2018interacting}). This is clearly the minimal sparsity threshold for which we can expect a mean field behavior, as evidenced by recent work in the extremely sparse (diluted) regime where $np_n$ converges to a finite non-zero constant (\cite{lacker2019large,oliveira2019interacting}).

The approach of \cite{delarue2017mean} given for the \Erdos \ case, based on the master equation, shows promise for a more general theory for dense graph models. But there are many difficulties to overcome, particularly in obtaining optimal sparsity thresholds as in our case study. It seems that the arguments of \cite{delarue2017mean} may extend to any graph sequence (not necessarily transitive) satisfying the conclusion of Proposition 7 therein, which is a fairly strong denseness assumption shown so far to cover the \Erdos \ case only when $p_n=p$ is constant. We do not pursue this any further, and we note also that \cite{delarue2017mean} deals with open-loop equilibria, whereas we work here with closed-loop.

\subsubsection{Random regular graphs.}

Random regular graphs are well known to admit tractable large-scale behavior, in both the dense and sparse regimes. Let $d_n \in \N\setminus\{1\}$, and let $R(n,d_n)$ denote a uniformly random choice of $d_n$-regular graph on $n$ vertices. In the dense regime, where  $d_n \to \infty$, Theorem \ref{th:approxeq} lets us construct approximate equilibria.
In the sparse regime, when $d_n =d$ is constant, the empirical measure $\mu_{R(n,d)}$ is known to converge weakly to an explicit continuous probability measure $\mu(d\lambda)$ known as (an affine transformation of) the Kesten-McKay law (\cite{kesten1959symmetric,mckay1981expected}), with density given by
\begin{align*}
\lambda \mapsto \frac{\sqrt{4(d-1)-d^2(\lambda+1)^2}}{2\pi(1 - (\lambda+1)^2)} 1_{\{|1+\lambda| \le 2\sqrt{d-1}/d\}}.
\end{align*}
The same limit $\mu$ arises for any sequence $G_n$ of $d$-regular graphs satisfying $\lim_{n\to\infty}\tfrac{C_k(G_n)}{|G_n|} = 0$ for each $k \in \N$, where $C_k(G_n)$ is the number of cycles of length $k$ in $G_n$, by \cite{mckay1981expected}. Note that we cannot apply our main results to the random regular graph $G_n=R(n,d)$ for $d$ fixed, because $G_n$ is then transitive with vanishing probability as $n\to\infty$; in fact, $G_n$ has trivial automorphism group with high probability (\cite{kim2002asymmetry}).

\subsection{The cooperative game.} \label{se:cooperative}

\begin{figure}[h]
    \centering
    \includegraphics[width=0.7\columnwidth]{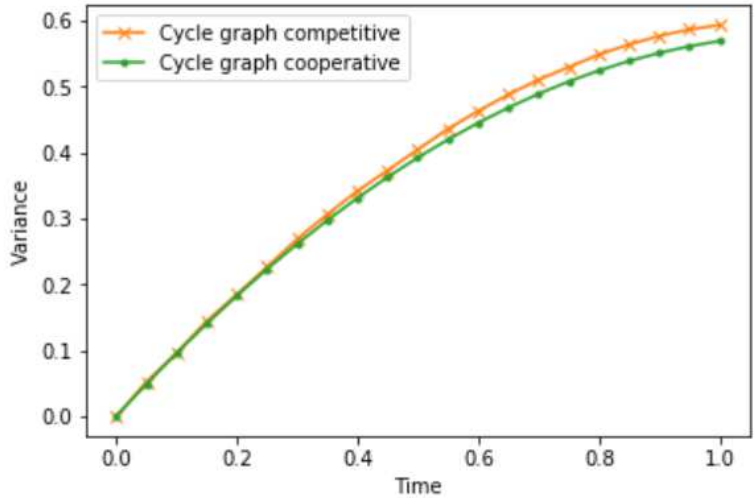}
    \caption{Variance of a typical player over time in the cycle graph, in the competitive versus cooperative regimes. Here $T=\sigma=c=1$.}
    \label{fig:compVScoop}
\end{figure}

For comparison, we discuss the corresponding cooperative game, which can be solved easily even without assuming transitivity of the underlying finite graph $G=(V,E)$.  In the setup of Section \ref{se:finitehorizon}, consider the optimal control problem
\begin{align*}
\inf_{\bm\alpha\in\A_G^V} \sum_{v \in V}J_v^G(\bm\alpha). 
\end{align*}
Let us abbreviate $\lap=\lap_G$ for the Laplacian.
The corresponding Hamilton-Jacobi-Bellman (HJB) equation is
\begin{equation*}
\left\{
\begin{array}{ll}
     &\partial_t v (t,\bm x) - \frac{1}{2} |\nabla v(t,\bm x)|^2 + \frac12 \sigma^2 \Delta v(t,\bm x) = 0, \quad (t,\bm x) \in (0,T) \times \R^V,  \\
     & v(T,\bm x) = \frac{c}{2}|\lap \bm x|^2 = \frac{c}{2}\bm x^\top \lap^\top \lap \bm x,
\end{array}
\right.   
\end{equation*}
and the optimal control is $\bm{\alpha}^* = - \nabla v$. Using the ansatz $v(t,\bm x) = \tfrac12 \bm x^\top F(t) \bm x + h(t)$,
for some symmetric matrix $F(t)$, 

the HJB equation becomes
\begin{equation*}
\frac{1}{2} \bm x^T F'(t)\bm x + h'(t) - \frac{1}{2} \bm x^\top F^2(t) \bm x + \frac12\sigma^2 \tr(F(t)) = 0, \quad (t,\bm x) \in (0,T) \times \R^V,
\end{equation*}
with terminal conditions $F(T) = c\lap^\top \lap$ and $h(T) = 0$. Matching coefficients, we deduce that $F$ and $h$ must solve
\begin{equation*}
    \left\{ \begin{array}{ll}
         & F'(t) - F^2(t) = 0, \hspace{0.5cm} F(T) = c\lap^\top \lap, \\
         & h'(t) + \frac12 \sigma^2 \tr(F(t)) = 0, \hspace{0.5cm} h(T) = 0.
    \end{array}
    \right.
\end{equation*}
We find that the solution to this system is given by
\begin{align*}
F(t) &= c \lap^\top \lap (I  + c(T-t) \lap^\top \lap)^{-1} = - \frac{d}{dt} \log(I+c(T-t)\lap^\top\lap),
\end{align*}
where the log of the positive definite matrix is defined via power series, and
\begin{align*}
h(t) &= \frac{\sigma^2}{2} \tr \log (I + c(T-t)\lap^\top \lap ) = \frac{\sigma^2}{2} \log  \det (I + c(T-t)\lap^\top \lap ).
\end{align*}
The optimal control is $\bm\alpha(t,\bm x)=-F(t)\bm x$, and the optimal state process follows 
\begin{equation*}
d\bm{X}(t) = - F(t) \bm{X}(t) + \sigma d\bm{W}(t).
\end{equation*}
This SDE can be explicitly solved, and the law of $\bm X(t)$ is Gaussian with mean $(I+cT\lap^\top\lap)(I+c(T-t)\lap^\top\lap)^{-1}\bm{X}(0)$ and covariance matrix
\begin{align*}
\sigma^2 \int_0^t (I+c(T-t)\lap^\top\lap)^2(I+c(T-s)\lap^\top\lap)^{-2} ds.
\end{align*}
In particular, if the graph $G$ is transitive, we compute for each $i \in V$ as in the beginning of Section \ref{se:asymptotics} that
\begin{align*}
\Var(X_i(t)) = \sigma^2 \int_0^t\int_{[-2,0]}\left(\frac{1+c(T-t)\lambda^2}{1+c(T-s)\lambda^2}\right)^2\mu_G(d\lambda)ds.
\end{align*}
And if $\bm X(0)=0$, then the per-player value is
\begin{align*}
\frac{1}{|V|}\inf_{\bm\alpha\in\A_G^V} \sum_{v \in V}J_v^G(\bm\alpha) = \frac{1}{|V|}v(0,0) = \frac{\sigma^2}{2} \int_{[-2,0]} \log (1 + cT\lambda^2) \,\mu_G(d\lambda).
\end{align*}
It is interesting to compare these outcomes to the competitive equilibrium computed in Theorem  \ref{thm:convergence_X^Gn_k}. In the competitive case, the function $f(t)=f_\mu(t)$ may be seen as determining the \emph{rate of flocking}; as $f(t)$ increases over time, players expend more effort to move toward the average. In the competitive case, this function depends crucially on the graph, and it simplifies to $f(t)=ct$ in the dense graph case. In the cooperative case, we always have $f(t)=ct$, but the solution is governed by the \emph{squared} Laplacian instead of the Laplacian itself. 
Figure \ref{fig:compVScoop} shows the variance over time for the competitive and cooperative solutions on the (limiting) cycle graph.

In addition to being an interesting point of comparison for the competitive case, we highlight the cooperative model also in connection with \emph{mean field control theory}. We are not the first to study stochastic optimal control problems on large graphs, but we do not know of much other work other than the recent papers of \cite{gao2017control,gao2018graphon}, which focus on the graphon setting and do not seem to offer any explicit examples.

\subsection{Comments on intrinsic limit models.} \label{se:intrinsic}

A key strength of the MFG paradigm is that it offers an intrinsic model, which gives information about the $n$-player equilibrium for large $n$ but which can be analyzed without reference to the latter.
There is no such intrinsic limiting model yet in the graph-based setting, aside from the recently proposed graphon MFGs for dense regimes (\cite{caines2018graphon,gao2020lqg}).
In the dense regime, our Theorem \ref{th:approxeq} says that the usual MFG model already suffices for the purpose of constructing approximate equilibria.
The sparse regime is more mysterious, and this section comments on some possibilities and difficulties.

In the sparse regime, it is natural to understand large graphs using the well-developed notion of \emph{local weak convergence}. See \cite{van2009randomII} for a thorough treatment, as well as \cite{lacker2019large,oliveira2019interacting}  for some recent work applying this framework to analyze large-scale interacting diffusion models on sparse graphs.
We will not define local weak convergence of graphs $G_n \to G$ here, but we highlight that it is known to imply the weak convergence of $\mu_{G_n}$ to a certain \emph{spectral measure} $\mu_G$  (\cite{BordenaveLelarge}). This spectral measure is defined intrinsically on any finite or countable locally finite graph $G=(V,E)$, by applying the spectral theorem to an appropriately defined Laplacian operator on the complex Hilbert space $\ell^2(V)$. In the finite (transitive) graph case, it coincides exactly with the empirical eigenvalue distribution defined in \eqref{def:spectralmeasure}.

Moreover, this formalism suggests an intrinsic description of the $n\to\infty$ limit, in the context of  Theorem \ref{thm:convergence_X^Gn_k}, in the sparse case. If $G_n$ converges in the local weak sense to an infinite graph $G$, we should expect the limit given in Theorem \ref{thm:convergence_X^Gn_k} to coincide with the equilibrium of the game set on $G$, defined as in Section \ref{se:finitehorizon} but with some extra care to handle the infinite population.
The equilibrium solution given in Theorem \ref{thm_equilibrium_n_player} should remain valid in the infinite transitive graph case, using the Laplacian \emph{operator} in place of the Laplacian \emph{matrix}, 
and being careful to interpret the countably infinite controlled SDE system in an appropriate (mild) sense.

But it is not at all clear what we stand to gain from this abstraction.
The usual MFG model, while often described as a model of a continuum of agents, is really characterized in terms of a \emph{single} representative player. The natural candidate for a sparse graph limit model, on the other hand, appears to be a game with infinitely many players, and it is not clear if any simpler, low-dimensional characterization is possible, even in the transitive case. The recent paper of \cite{lacker2019large} derives low-dimensional marginal dynamics for uncontrolled models on regular (and unimodular Galton-Watson) trees, but this is only possible when the dynamics are \emph{local} in the sense that the drift of $X^G_v$ depends only on the neighbors $(X^G_u)_{u \sim v}$. This is not the case in equilibrium in our model. 

The issue of non-local equilibrium controls suggests an alternative strategy of truncating the range of dependence, with each player basing its control only on those other players within some pre-specified graph distance. This may give rise to a sequence of more tractable approximate equilibria. We leave this and the other speculations of this section for future work.

{\ }

\noindent\textbf{Organization of the rest of the paper.} The rest of the paper is devoted to the proofs of the results stated in this section. We begin in Section \ref{se:ODEanalysis} with an analysis of the functions $Q_G$ and $Q_\mu$ and the ODEs for $f_G$ and $f_\mu$ which appear in Theorems \ref{thm_equilibrium_n_player} and \ref{thm:convergence_X^Gn_k}. Section \ref{section:n_player_finite} then gives the proof of Theorem \ref{thm_equilibrium_n_player}. This proof is essentially a verification argument and does not explain how to derive the announced solution, so in Section \ref{se:alternateproof} we give a sketch of a direct derivation. Section \ref{subsection:cov_bounds} proves the covariance bound of Proposition \ref{pr:corrdecay}. Finally, Section \ref{se:asymptoticregimes} proves Theorem \ref{thm:convergence_X^Gn_k}, Theorem \ref{th:approxeq}, and some claims of Section \ref{se:examples}.

\section{Analysis of the ODE} \label{se:ODEanalysis}

In this section, we derive several useful results about the ODEs encountered in Theorems \ref{thm_equilibrium_n_player} and \ref{thm:convergence_X^Gn_k}. For a probability measure $\mu$ on $[-2,0]$, define the function
\begin{align}
Q_\mu(x) = \exp\int_{[-2,0]} \log(1-x\lambda)\,\mu(d\lambda), \quad x \ge 0. \label{def:Qmu}
\end{align}
Note that the support of $\mu$ ensures that $Q_\mu(x)$ is well-defined and infinitely differentiable for $x \ge 0$. Note that if $\mu=\mu_G$ for a finite graph $G$, recalling the notation of Section \ref{se:asymptotics}, then $Q_\mu=Q_G$ takes the form of a normalized determinant as in \eqref{def:intro:Q}. 

We restrict to probability measures on $[-2,0]$ in this section precisely because $\mu_G$ is supported on $[-2,0]$ for every finite graph $G$, as discussed in Remark \ref{re:spectrum}. In addition, because the adjacency matrix of a (simple) graph has zero trace, the Laplacian matrix of a graph on $n$ vertices therefore has trace $-n$. In particular,
\begin{align}
\int_{[-2,0]} x\,\mu_G(dx) = \frac{1}{n}\sum_{i=1}^n\lambda^G_i = \frac{1}{n}\tr(\lap_G)=-1,
\end{align}
for a finite graph $G$ on $n$ vertices with Laplacian eigenvalues $(\lambda^G_1,\ldots,\lambda^G_n)$.
We thus restrict our attention to the set $\P_{\mathrm{Lap}}$ of probability measures on $[-2,0]$ with mean $-1$, equipped with the topology of weak convergence (i.e., $\mu_n \to \mu$ in $\P_{\mathrm{Lap}}$ if $\int f\,d\mu_n \to \int f\,d\mu$ for each continuous function $f : [-2,0] \to \R$). 

We will make repeated use of the following formulas for the first two derivatives of $Q_\mu$, computed via straightforward calculus:
\begin{align}
Q'_\mu(x) &= Q_\mu(x)\int_{[-2,0]} \frac{-\lambda}{1-x\lambda}\,\mu(d\lambda), \label{pf:Qderiv} \\
Q''_\mu(x) &= Q_\mu(x)\left[\left(\int_{[-2,0]} \frac{-\lambda}{1-x\lambda}\,\mu(d\lambda)\right)^2 - \int_{[-2,0]} \frac{\lambda^2}{(1-x\lambda)^2}\,\mu(d\lambda) \right]. \nonumber
\end{align}

\begin{proposition} \label{pr:Qproperties}
For each $\mu \in \P_{\mathrm{Lap}}$,
the function $Q_\mu : \R_+ \to \R_+$ defined in \eqref{def:Qmu} satisfies
\begin{align*}
1 &\le Q_\mu(x) \le 1+x, \qquad 0 < Q'_\mu(x) \le 1, \qquad -4 \le Q''_\mu(x) \le 0, 
\end{align*}
for all $x \in \R_+$, as well as
\begin{align}
Q'_\mu(x) &\ge 1 - (x + \tfrac12 x^2)\Var(\mu), \quad \text{where}\quad  \Var(\mu) := \int_{[-2,0]} (\lambda+1)^2\mu(d\lambda). \label{def:QsharpLB}
\end{align}
\end{proposition}

\proof{Proof of Proposition \ref{pr:Qproperties}.}
First note that the support of $\mu$ ensures that $Q_\mu(x) \ge 0$ for all $x \ge 0$.
Jensen's inequality yields
\begin{align*}
Q_\mu(x) \le \int_{[-2,0]} (1-x\lambda)\,\mu(d\lambda) = 1 + x.
\end{align*}
and $Q''_\mu(x) \le 0$. Since $Q'_\mu(0)=1$, we deduce that $Q'_\mu(x) \le 1$ for all $x \ge 0$.
The next claims follow from the fact that, for each $x \ge 0$, the function
\begin{align}
[-2,0] \ni \lambda \mapsto -\lambda/(1-x\lambda) \quad \text{is nonnegative and strictly decreasing.} \label{pf:gdecreasing}
\end{align}
Indeed, this first implies that $Q'_\mu(x) \ge 0$, and in fact the inequality must be strict because $\mu$ has mean $-1$ and is thus not equal to $\delta_0$.
In addition, \eqref{pf:gdecreasing} implies that 
\begin{align*}
Q_\mu''(x) &\ge -Q_\mu(x)\int_{[-2,0]} \frac{\lambda^2}{(1-x\lambda)^2}\,\mu(d\lambda) \ge -Q_\mu(x) \frac{4}{(1+2x)^2} \ge -4,
\end{align*}
where the last step uses $Q_\mu(x) \le 1+x\le (1+2x)^2$.
To prove the final claim, define
\[
\theta(x):=\int_{[-2,0]}\frac{-\lambda}{1-x \lambda}\,\mu(d\lambda), \quad x \ge 0.
\]
Using the Harris inequality (known as Chebyshev's sum inequality in the discrete case), we get
\begin{align*}
-2\theta(x)\theta'(x) &= \int_{[-2,0]}\frac{-\lambda}{1-x \lambda}\,\mu(d\lambda)\int_{[-2,0]}\frac{2\lambda^2}{(1-x \lambda)^2}\,\mu(d\lambda) \\
	&\le \int_{[-2,0]}\frac{-2\lambda^3}{(1-x \lambda)^3}\,\mu(d\lambda) = \theta''(x),
\end{align*}
since both integrands are decreasing functions of $\lambda$. Note that $\theta(0)=1$ and $\theta'(0)=\Var(\mu)+1$, and integrate the above inequality to find $\theta^2(x) + \theta'(x) \ge -\Var(\mu)$ for $x \ge 0$. The identity $Q''_\mu(x)=Q_\mu(x)(\theta^2(x)+\theta'(x))$ thus yields $Q''_\mu(x) \ge -(1+x)\Var(\mu)$. Integrate, using $Q'_\mu(0)=1$, to complete the proof.
\endproof

Using these properties, we next justify the existence, uniqueness, and stability for the ODEs appearing in Theorems \ref{thm_equilibrium_n_player} and \ref{thm:convergence_X^Gn_k}:

\begin{proposition} \label{pr:f-convergence}
Let $\mu \in \P_{\mathrm{Lap}}$. There is a unique continuous function $f_\mu : \R_+ \to \R_+$, continuously differentiable on $(0,\infty)$, satisfying
\begin{align*}
f'_\mu(t) = cQ'_\mu(f_\mu(t)), \ \ \ t > 0, \qquad f_\mu(0)=0.
\end{align*}
Moreover, we have the bounds
\begin{align*}
0 \le f_\mu(t) \le ct, \qquad f_\mu(t) \ge ct - (\tfrac12 c^2t^2 + \tfrac16 c^3t^3)\Var(\mu).
\end{align*}
Finally, if $\mu_n$ is a sequence in $\P_{\mathrm{Lap}}$ converging to $\mu$, then $Q'_{\mu_n}$ and $f_{\mu_n}$ converge uniformly to $Q'_{\mu}$ and $f_{\mu}$, respectively, on compact subsets of $\R_+$.
\end{proposition}

\proof{Proof of Proposition \ref{pr:f-convergence}.}
By Proposition \ref{pr:Qproperties}, $Q'_\mu$ is nonnegative and Lipschitz (with constant $4$) on $\R_+$. A standard Picard iteration  yields existence and uniqueness of the ODE in question.

Next we prove the estimates for $f_\mu$. Note first that the bound $Q'_\mu \le 1$ from Proposition \ref{pr:Qproperties} ensures that $0 \le f_\mu(t) \le ct$ for all $t \ge 0$. 
Use the lower bound on $Q'_\mu$ from \eqref{def:QsharpLB} to get
\begin{align*}
f'_\mu(t) &= cQ'_\mu(f_\mu(t)) \ge c - c\left(f_\mu(t) + \frac12 f_\mu(t)^2\right) \Var(\mu)  \ge c - c(ct + \tfrac12 c^2t^2)\Var(\mu).
\end{align*}
Integrate both sides to get the desired lower bound on $f_\mu(t)$.

We next prove that $Q'_{\mu_n}$ converges to $Q'_{\mu}$ uniformly on compacts. Because $Q'_{\nu}$ is $4$-Lipschitz and $Q'_{\nu}(0)=1$ for each $\nu \in \P_{\mathrm{Lap}}$, the equicontinuous family $\{Q'_{\nu} : \nu \in \P_{\mathrm{Lap}}\} \subset C(\R_+;\R)$ is precompact in the topology of uniform convergence on compacts, by the Arzel\`a-Ascoli theorem. Hence, we need only prove the pointwise convergence of $Q'_{\mu_n}$ to $Q'_{\mu}$. But this follows easily from the assumed convergence $\mu_n \to \mu$ and the form of $Q'$ in \eqref{pf:Qderiv}.
Finally, since $Q'_{\mu_n}$ is $4$-Lipschitz, for $t \ge 0$ we have
\begin{align*}
|f_{\mu_n}(t)-f_\mu(t)| &\le c\int_0^t|Q'_{\mu_n}(f_{\mu_n}(s)) - Q'_{\mu_n}(f_{\mu}(s))|ds + c\int_0^t|Q'_{\mu_n}(f_{\mu}(s)) - Q'_{\mu}(f_\mu(s))|ds \\
	&\le 4c\int_0^t| f_{\mu_n}(s) - f_\mu(s) |ds + c\sup_{u \in [0,ct]}|Q'_{\mu_n}(u) - Q'_{\mu}(u)|.
\end{align*}
Use Gronwall's inequality and the uniform convergence of $Q'_{\mu_n}$ to $Q'_{\mu}$ on compacts to deduce that $f_{\mu_n} \to f_\mu$ uniformly on compacts.
\endproof

This concludes the basic analysis of $f_\mu$ and $Q_\mu$ needed for the proofs of the main results of Sections \ref{se:results_nplayer_finite} and \ref{se:asymptotics}. The rest of the section is devoted to some properties of the variance computed in Theorem \ref{thm:convergence_X^Gn_k}, which will not be needed in the subsequent sections but which justify the claims in Remarks \ref{re:variance} and \ref{re:variance-dense}.
Extending the formula \eqref{def:gaussianlimit} for  the variance, we define for $\mu \in \P_{\mathrm{Lap}}$ and $t \in [0,T]$ the quantity
\begin{align}
\varfunc_\mu(t) := \sigma^2\int_0^t\int_{[-2,0]} \left(\frac{1 - \lambda f_\mu(T-t)}{1 - \lambda f_\mu(T-s)}\right)^2\mu(d\lambda) ds. \label{def:vmu}
\end{align}

\begin{proposition} \label{pr:fQproperties}
For each $\mu \in \P_{\mathrm{Lap}}$, $f_\mu$ is strictly increasing and concave, with $0 < f'_\mu(t) \le c$ for all $t \ge 0$.  Moreover, $\varfunc_\mu$ satisfies
\begin{align*}
\varfunc_\mu(0)=0, \quad \varfunc'_\mu(0)=\sigma^2, \quad \text{and } \ \varfunc''_\mu(0) = -2\sigma^2\frac{(f'_\mu(T))^2}{cQ_\mu(f_\mu(T))}.
\end{align*}
\end{proposition}
\proof{Proof of Proposition \ref{pr:fQproperties}.}
Fix $\mu \in \P_{\mathrm{Lap}}$. First, we note that
\begin{align*}
f'_\mu(t)=cQ'_\mu(f_\mu(t)), \qquad \text{and} \qquad f''_\mu(t) &= cQ''_\mu(f_\mu(t))f'_\mu(t).
\end{align*}
Using Proposition \ref{pr:Qproperties} we deduce that $ 0 <f'_\mu(t) \le c$ and $f''_\mu(t) \le 0$ for all $t \ge 0$, which proves the first claim. 
Differentiating in \eqref{def:vmu} using Liebniz's rule yields
\begin{align*}
\varfunc'_\mu(t) &= \sigma^2 - 2\sigma^2f'_\mu(T-t) \int_0^t\int_{[-2,0]} \frac{-\lambda(1 - \lambda f_\mu(T-t))}{(1 - \lambda f_\mu(T-s))^2}\mu(d\lambda) ds.
\end{align*}
The claim $\varfunc'_\mu(0)=\sigma^2$ follows. Differentiate again at $t=0$ to get
\begin{align*}
\varfunc''_\mu(0) &= - 2\sigma^2f'_\mu(T)  \int_{[-2,0]} \frac{-\lambda }{ 1 - \lambda f_\mu(T) }\mu(d\lambda) \\
	&= - 2\sigma^2\frac{d}{dt}\Big|_{t=T} \int_{[-2,0]} \log( 1 - \lambda f_\mu(t))\mu(d\lambda) \\
	&= - 2\sigma^2\frac{d}{dt}\Big|_{t=T} \log Q_\mu(f_\mu(t)) \\
	&= - 2\sigma^2 \frac{Q'_\mu(f_\mu(T))f'_\mu(T)}{Q_\mu(f_\mu(T))} \\
	&= - 2\sigma^2 \frac{(f'_\mu(T))^2}{cQ_\mu(f_\mu(T))}.
\end{align*}
\endproof

We next show that the dense graph regime uniquely minimizes the variance, as announced in Remark \ref{re:variance-dense}.

\begin{proposition} \label{pr:vardeltamin}
Let $\mu \in \P_{\mathrm{Lap}}$. For each $t \in (0,T]$ we have
\begin{align*}
\varfunc_\mu(t) \ge \sigma^2t \frac{1+c(T-t) }{1+cT} = \varfunc_{\delta_{-1}}(t),
\end{align*}
with equality if and only if $\mu=\delta_{-1}$.
\end{proposition}
\proof{Proof of Proposition \ref{pr:vardeltamin}.}
Note from Proposition \ref{pr:Qproperties} that $Q_\mu(x) \le 1+x = Q_{\delta_{-1}}(x)$ and $Q'_\mu(x) \le 1 = Q'_{\delta_{-1}}(x)$  for $x \ge 0$. By a standard comparison argument, it follows that $f_\mu(t) \le f_{\delta_{-1}}(t)=ct$ for all $t \ge 0$. Next, note that
\begin{align}
\R_+ \ni t \mapsto \frac{1-\lambda f_\mu(t)}{1-\lambda f_{\delta_{-1}}(t)} = \frac{1-\lambda f_\mu(t)}{1-\lambda ct} \ \ \text{ is non-increasing for each } \lambda \in [-2,0]. \label{pf:nondecreasing1}
\end{align}
Indeed, the derivative is
\begin{align*}
\frac{ -\lambda f'_\mu(t)(1-\lambda ct) + \lambda c(1-\lambda f_\mu(t))}{(1-\lambda ct)^2} = \frac{-\lambda (f'_\mu(t) - c) + c\lambda^2(tf'_\mu(t) - f_\mu(t)) }{(1-\lambda ct)^2}.
\end{align*}
This is at most zero, because we know $f'_\mu(t) \le c$ by Proposition \ref{pr:fQproperties}, and the concavity of $f_\mu$ together with $f_\mu(0) = 0$ imply $tf'_\mu(t) \le f_\mu(t)$. This proves \eqref{pf:nondecreasing1}, which is equivalent to the fact that
\begin{align*}
\frac{1-\lambda f_\mu(T-t)}{1-\lambda f_\mu(T-s)} \ge \frac{1-\lambda f_{\delta_{-1}}(T-t)}{1-\lambda f_{\delta_{-1}}(T-s)} = \frac{1-\lambda c(T-t)}{1-\lambda c(T-s)}, \quad \text{for } t > s > 0.
\end{align*}
As both sides are non-negative, this implies
\begin{align}
\int_{[-2,0]} \left(\frac{1 - \lambda f_\mu(T-t)}{1 - \lambda f_\mu(T-s)}\right)^2\mu(d\lambda) \ge \int_{[-2,0]} \left(\frac{1-\lambda c(T-t)}{1-\lambda c(T-s)} \right)^2\mu(d\lambda). \label{pf:vmin1}
\end{align}
Finally, the function
\begin{align*}
[-2,0]\ni \lambda \mapsto \left(\frac{ 1- \lambda c(T-t)}{1- \lambda c(T-s)}\right)^2
\end{align*}
is strictly convex for $t > s > 0$, which by Jensen's inequality implies 
\begin{align}
\int_{[-2,0]} \left(\frac{1-\lambda c(T-t)}{1-\lambda c(T-s)} \right)^2\mu(d\lambda) \ge \left(\frac{1 +c(T-t)}{1 +c(T-s)}\right)^2 \label{pf:vmin2}
\end{align}
since $\mu$ has mean $-1$. Combine \eqref{pf:vmin1} and \eqref{pf:vmin2} with the definition of $\varfunc_\mu$ to get $\varfunc_\mu(t) \ge \varfunc_{\delta_{-1}}(t)$, as desired. The inequality \eqref{pf:vmin2} is strict unless $\mu=\delta_{-1}$.
\endproof

\section{The equilibrium on finite graphs: Proof of Theorem \ref{thm_equilibrium_n_player}}\label{section:n_player_finite}

This section gives the proof of Theorem \ref{thm_equilibrium_n_player}. 
We begin with some general symmetry considerations in Section \ref{se:symmetries}. We then derive the HJB system in Section \ref{se:HJBs} and reduce it to a system of Riccati equations in Section \ref{se:Riccati}; these two steps are standard for linear-quadratic games. The explicit resolution of the system of Riccati equations is where the difficulty lies. In Section \ref{se:checkingRiccati} we show that the proposed solution Theorem \ref{thm_equilibrium_n_player} does indeed work, and the remaining sections \ref{se:stateprocess} and \ref{se:valuecomputation} provide the remaining computations of the equilibrium state process dynamics and the average value of the game.

The proof given in this section, while complete and rigorous, is opaque in the sense that it does not give any idea of how one might arrive at the solution of the system of Riccati equations. For this reason, we give in Section \ref{se:alternateproof} a sketch a direct derivation of the solution.

We fix throughout the section a finite transitive graph $G = (V,E)$, and write $V=\{1,\ldots,n\}$ for some $n \in \N$. We may assume without loss of generality that $G$ is connected (see Remark \ref{re:disconnected}). 
Throughout this entire section, we omit $G$ from the notation by writing, e.g., $\lap=\lap_G$ and $\bm X=\bm X^G$. 

\subsection{Symmetries.} \label{se:symmetries}

We first discuss some basic symmetry properties. Since the graph $G$ is transitive and thus regular, the Laplacian matrix is symmetric, i.e., $\lap=\lap^\top$. 

We will make some use of the so-called \emph{regular representation} of the automorphism group. Recall from the beginning of Section \ref{se:results_nplayer_finite} that $\Aut(G)$ denotes the set of automorphisms of $G$. To each $\varphi \in \Aut(G)$ we associate an invertible $n \times n$ matrix $R_\varphi$, defined by requiring $R_\varphi e_i = e_{\varphi(i)}$ for each $i \in V$, where we recall that $(e_1, ..., e_n)$ denotes the standard Euclidean basis in $\R^n$. It is clear that $R_\varphi R_\psi = R_{\varphi \circ \psi}$, and in particular $R_\varphi^{-1} = R_{\varphi^{-1}}$. We also have $R_\varphi^\top = R_{\varphi^{-1}}$, because
\begin{align*}
e_i^\top R_\varphi e_j &= e_i^\top e_{\varphi(j)} = 1_{\{i = \varphi(j)\}} = 1_{\{j = \varphi^{-1}(i)\}} = e_j^\top R_{\varphi^{-1}}e_i.
\end{align*}
The following elementary lemma summarizes some uses of transitivity (Definition \ref{def:transitive}). The third property will be used only in the alternative proof of Theorem \ref{thm_equilibrium_n_player} given in Section \ref{se:alternateproof}.

\begin{lemma} \label{le:symmetries}
Assume that $G=(V,E)$ is transitive, with $V=\{1,\ldots,n\}$.
\renewcommand\theenumiii{\@roman\c@enumiii}
\begin{enumerate}[(i)]

\item $\lap$ commutes with $R_\varphi$ for each $\varphi \in \Aut(G)$.
\item If $Y \in \R^{n \times n}$ commutes with $R_{\varphi}$ for every $\varphi \in \Aut(G)$, then
\begin{equation*}
Y_{ii} = \frac{1}{n}\tr(Y), \qquad \forall \, i\in V.
\end{equation*}
\item If $Y^1, \ldots, Y^n \in \R^{n \times n}$ satisfy $R_{\varphi} Y^i = Y^{\varphi(i)}R_{\varphi}$ for every $\varphi \in \Aut(G)$ and $i \in V$, then
\begin{equation*}
Y^{i}_{ii} = Y^{j}_{jj}, \qquad  \forall \, i,j \in V.
\end{equation*}
\end{enumerate}
\end{lemma}
\proof{Proof of Lemma \ref{le:symmetries}.} {\ }
\begin{enumerate}[(i)]
\item For $\varphi \in \Aut(G)$ we compute
\begin{align*}
e_i^\top R_\varphi^\top \lap R_\varphi e_j = e_{\varphi(i)}^\top \lap e_{\varphi(j)} = e_i^\top \lap e_j, \qquad \forall i,j \in V,
\end{align*}
with the last equality using the fact that $\varphi$ is an automorphism. This shows that $R_\varphi^\top \lap R_\varphi=\lap$. Since $R_\varphi^\top=R_\varphi^{-1}$, this shows $\lap R_{\varphi} = R_{\varphi}\lap$.
\item It suffices to apply (iii) with $Y^1=\cdots=Y^n=Y$ to get $Y_{ii}=Y_{jj}$ for all $i,j \in V$.
\item Let $i,j \in V$. By transitivity, there exists $\varphi \in \Aut(G)$ such that $\varphi(i) = j$. Then
\begin{equation*}
e_i^\top Y^i e_i = e_i^\top R_{\varphi}^T Y^{\varphi(i)} R_{\varphi} e_i = e_{\varphi(i)}^\top Y^{\varphi(i)} e_{\varphi(i)} = e_j^\top Y^j e_j.
\end{equation*}
\end{enumerate}
\endproof

\subsection{The corresponding system of HJB equations.}\label{se:HJBs}

We can write the cost function of \eqref{def:cost} for player $i \in V$ as
\begin{align*}
J_i(\alpha_1,\ldots,\alpha_n)=\frac12\E\left[ \int_0^T |\alpha_i(t, \bm{X}(t))|^2dt + c\left|e_i^\top\lap\bm X(T)\right|^2 \right].
\end{align*}
A standard argument associates this $n$-player game to the following Nash system of $n$ coupled PDEs:
\begin{align}
\begin{split}
&0 = \partial_t v_i(t,\bm x) - \frac12(\partial_i v_i(t,\bm x))^2 - \sum_{k \ne i} \partial_k v_k(t,\bm x) \partial_k v_i(t,\bm x) + \frac{\sigma^2}{2} \sum_{k=1}^n \partial_{kk} v_i(t,\bm x),  \\
&v_i(T,\bm x) = \frac12 c(e_i^\top\lap \bm x)^2, \qquad i=1,\ldots,n.
\end{split} \label{n-player-HJB}
\end{align}
Here, we write $\bm x=(x_1,\ldots,x_n)$ for a typical vector in $\R^n$, and for the functions $v_i : [0,T] \times \R^n \to \R$ we write $\partial_t$ and $\partial_k$ for the derivative with respect to $t$ and $x_k$, respectively.
If $(v_1,\ldots,v_n)$ is a classical solution of \eqref{n-player-HJB}, then the controls
\begin{align}
\alpha_i(t,\bm x) = -\partial_i v_i(t,\bm x), \quad i =1,\ldots,n, \label{n-HJB-alphas}
\end{align}
form a Markovian Nash equilibrium.

For a thorough derivation of PDEs of this form and a verification theorem, we refer to the book of \cite[Section I.2.1.4]{carmona-delarue-book}.
But let us briefly explain how to check if some $(\alpha_1,\ldots,\alpha_n) \in \A_G^n$ forms a Nash equilibrium. Considering player $i$'s optimization problem,
\begin{align}
\inf_{\alpha_i \in \A_G}J_i(\alpha_1,\ldots,\alpha_{i-1},\alpha_i,\alpha_{i+1},\ldots,\alpha_n). \label{oneplayerOPT}
\end{align}
Standard stochastic control theory (see, e.g., \cite{fleming2006controlled,pham2009continuous}) leads to the HJB equation
\begin{align}
\begin{split}
&0 = \partial_t v_i(t,\bm x) - \frac12(\partial_i v_i(t,\bm x))^2 + \sum_{k \ne i} \alpha_k(t,\bm x) \partial_k v_i(t,\bm x) + \frac{\sigma^2}{2} \sum_{k=1}^n \partial_{kk} v_i(t,\bm x),  \\
&v_i(T,\bm x) = \frac12 c(e_i^\top\lap \bm x)^2.
\end{split} \label{oneplayerHJB}
\end{align}
Indeed, the (reduced) Hamiltonian for player $i$ is
\begin{align*}
H_i(\bm x,\bm p) = \inf_{a \in \R}\left(\frac12 a^2 + a p_i + \sum_{k \neq i}\alpha_k(t,\bm x)p_k\right) =  - \frac{1}{2} p_i^2 + \sum_{k \neq i}\alpha_k(t,\bm x)p_k, \qquad \bm x, \bm p \in \R^n,
\end{align*}
with the infimum attained at $a=-p_i$. Hence, after solving the PDE \eqref{oneplayerHJB}, the optimal control in \eqref{oneplayerOPT} is given by $\alpha_i(t, \bm x) = -\partial_iv_i(t,\bm x)$. Applying this optimality criterion for each player $i=1,\ldots,n$ couples the PDEs \eqref{oneplayerHJB}, leading to the system \eqref{n-player-HJB}.

\subsection{Reduction to Riccati equations.} \label{se:Riccati}

Linear-quadratic control problems and games always relate to Riccati-type equations after a quadratic ansatz. For our PDE system \eqref{n-player-HJB}, we make the ansatz
\begin{align}
v_i(t,\bm x) = \frac12 \bm x^\top F^i(t) \bm x + h_i(t), \qquad (t,\bm x) \in [0,T] \times \R^n, \label{ansatz}
\end{align}
where $F^i : [0,T] \to \R^{n \times n}$ and $h_i : [0,T] \to \R$ are functions to be determined.
We assume without loss of generality that $F^i(t)$ is symmetric for each $t$.

\begin{lemma}\label{lem:Fieq}
Any solution of the Nash system \eqref{n-player-HJB} of the form \eqref{ansatz}, must satisfy the equations
\begin{align}
0 &= \Dot{F}^i(t)  -  \sum_{j = 1}^n F^j(t)e_j e_j^\top F^i(t) - F^i(t)\sum_{j = 1}^n e_j e_j^\top F^j(t) + F^i(t)e_i e_i^\top  F^i(t), \label{equation-Fi} \\
0 &= \Dot{h}_i(t) + \frac{\sigma^2}{2} \tr(F^i(t)), \qquad i=1,\ldots,n, \nonumber
\end{align}
and the boundary conditions
\begin{align}
F^i(T) = c\lap e_i e_i^\top \lap, \qquad h_i(T)=0. \label{ansatz:boundary}
\end{align}
\end{lemma}

\proof{Proof of Lemma \ref{lem:Fieq}.}
Recall that $\lap=\lap^\top$, and note that the boundary condition $v_i(T,\bm x)=\frac12 c(e_i^\top\lap \bm x)^2 = \frac{1}{2} c \bm x^\top \lap e_i e_i^\top \lap\bm x$ leads to the boundary conditions.
We write $\Dot{F^i}$ and $\Dot{h}_i$ for the derivatives of these functions. Once we check that this ansatz is correct, the equilibrium controls are given by
\begin{align}
\alpha_i(t,\bm x) =  - \partial_i v_i(t,\bm x)=  -  e_i^\top F^i(t)\bm x, \quad i=1,\ldots,n. \label{n-HJB-alphas-ansatz}
\end{align}

Applying the ansatz \eqref{ansatz} to the PDE system \eqref{n-player-HJB}, noting that $\partial_kv_i(t,\bm x) = e_k^\top F^i(t)\bm x$ and $\partial_{kk}v_i(t,\bm x) = e_k^\top F^i(t) e_k$, leads to
\begin{align*}
\frac{1}{2}\bm x^\top \Dot{F}^i(t) \bm x + \Dot{h}_i(t) + \frac{1}{2} (e_i^\top F^i(t) \bm x)^2 - \sum_{k = 1}^n \left(e_k^\top F^k(t) \bm x \right)\left(e_k^\top F^i(t) \bm x \right)  + \frac{\sigma^2}{2} \tr(F^i(t)) = 0.
\end{align*}
Collecting like terms, we find
\begin{align*}
\bm x^\top \Big( \frac12 \Dot{F}^i(t)  -  \sum_{k = 1}^n F^k(t)e_k e_k^\top F^i(t) +\frac12 F^i(t)e_i e_i^\top  F^i(t) \Big) \bm x + \Dot{h}_i(t) + \frac{\sigma^2}{2}\tr(F^i(t)) =0.
\end{align*}
This must hold for each $\bm x \in \R^n$, and we note that a square matrix $A$ satisfies $\bm x^\top A \bm x=0$ if any only if $A+A^\top=0$. We conclude by recalling that $F^i(t)$ is symmetric.

\endproof

Once we solve for $F^i$ using the first equation, the second equation and the boundary condition $h_i(T)=0$ yield
\begin{align}
h_i(t) &= \frac{\sigma^2}{2}\int_t^T \tr(F^i(s))ds. \label{solution-h_i}
\end{align}
Hence, the main task ahead is to solve the system \eqref{equation-Fi}.
A key role will be played by the matrix
\begin{align}
\widehat{P}(t) = \sum_{j = 1}^n F^j(t)e_j e_j^\top, \label{eq:Phat}
\end{align}
which appears multiplied by $F^i(t)$ in equation \eqref{equation-Fi}. In equilibrium, $\widehat{P}(t)$ 
will agree with $P_G(t)$ defined in \eqref{def:intro:P*}. If we freeze the function $\widehat{P}(t)$ in \eqref{equation-Fi}, the system of equations for $(F^1,\ldots,F^n)$ decouples, and each $F^i$ satisfies the simpler matrix Riccati equation
\begin{equation}\label{ODEsF*}
\Dot{F^i}(t) - \widehat{P}(t) F^i(t) - F^i(t) \widehat{P}(t) + \widehat{P}(t) e_i e_i^T \widehat{P}(t) = 0,  \qquad  t \in (0,T), \ \ i=1,\ldots,n,
\end{equation}
which we can solve explicitly in terms of $\widehat{P}(t)$. This direct strategy will be carried out in Section \ref{se:alternateproof}, ultimately leading to a fixed point equation \eqref{P_asfct_MPtilde} that $\widehat{P}(t)$ must satisfy. But we will first complete this section by showing that the solution proposed in Theorem \ref{thm_equilibrium_n_player} is indeed correct: We show that $P_G$, and suitable functions $(F^1,\ldots,F^n)$ thereof, do indeed simultaneously solve the equations \eqref{eq:Phat} and \eqref{ODEsF*}, thus providing a solution of the system in Lemma \ref{lem:Fieq}.

\subsection{Checking the solution.} \label{se:checkingRiccati}

It follows from the results of Section \ref{se:ODEanalysis} that the ODE
\begin{align*}
f'(t) = c Q'(f(t)), \qquad f(0) = 0,
\end{align*}
is well-posed, where $Q(x)=Q_G(x)=(\det(I-x\lap))^{1/n}$, and the solution $f$ is nonnegative and strictly increasing since $Q' > 0$.
As a result, the matrix-valued function defined in \eqref{def:intro:P*} by
\begin{align}
P(t) =P_G(t):= - f'(T-t) \lap \big(I - f(T-t) \lap \big)^{-1} \label{def:P}
\end{align}
is well-defined; the symmetric matrix $I - f(T-t) \lap$ is invertible because $\lap$ is negative semidefinite and $f \ge 0$.
Moreover, $P(t)$ satisfies the following useful properties for each $t \in [0,T]$:
\begin{enumerate}[(i)]
\item $e_i^T P(t) e_i = \frac{1}{n} \tr(P(t)) > 0$ for all $i = 1,\ldots, n$.
\item $P(t)$ is symmetric and positive semidefinite. Both $\lap$ and $P(t)$ have 0 as an eigenvalue, with the same multiplicity.
\item $P(t)$ and $\lap$ commute.
\end{enumerate}
Indeed, the third property is trivial. The second follows from the facts that $\lap$ is negative semidefinite, $f \ge 0$, and $f' > 0$; note that the vector of all ones is an eigenvector of $\lap$ with eigenvalue zero, and thus also an eigenvector of $P(t)$ with eigenvalue zero.
To derive the first property, note that $\lap$ commutes with $R_{\varphi}$ for all $\varphi \in \Aut(G)$ by Lemma \ref{le:symmetries}(i), and thus so does $P(t)$, which means Lemma \ref{le:symmetries}(ii) applies. The strict positivity of $\tr(P(t))$ follows from the fact that $P(t)$ is positive semidefinite and is not the zero matrix.

Since $\tr(P(t)) > 0$ by property (i), we may define $F^i$ for $i=1,\ldots,n$ by
\begin{align}
F^i(t) := \frac{1}{\tr(P(t))/n} P(t)e_ie_i^\top P(t). \label{pf:def:F*i}
\end{align}
Using property (i), we compute
\begin{align*}
e_i^\top F^i(t) &= \frac{1}{\tr{(P(t))}/n}\left(e_i^\top P(t) e_i\right)e_i^\top P(t) = e_i^\top P(t).
\end{align*}
In other words, the $i^\text{th}$ column of $F^i(t)$ is the same as that of $P(t)$. In particular, the control $\alpha^*_i(t,\bm x) = - e_i^\top P(t) \bm x$ defined in Theorem \ref{thm_equilibrium_n_player} satisfies $\alpha^*_i(t,\bm x) = - e_i^\top F^i(t) \bm x$. Recall from \eqref{n-HJB-alphas-ansatz} that this was indeed the form dictated by the PDE.

We next check that $(F^1,\ldots,F^n)$ solve the equations \eqref{equation-Fi} and boundary conditions \eqref{ansatz:boundary} identified in Lemma \ref{lem:Fieq}. 
Beginning with the boundary condition \eqref{ansatz:boundary}, note first that $f(0)=0$ and thus
\begin{equation*}
P(T) = - f'(0)\lap = - c Q'(0) \lap = - c \lap.
\end{equation*}
This implies $\frac{\tr{(P(T)})}{n} = c$, and from which we immediately compute $F^i(T)=c\lap e_ie_i^\top \lap$, as desired.
We next turn to the equations \eqref{equation-Fi}.
Using property (i) once again, we compute 
\begin{equation*}
\sum_{i=1}^n F^i(t)e_ie_i^\top = \frac{1}{\tr{(P(t))}/n} \sum_{i=1}^n P(t)e_i\left(e_i^\top P(t) e_i\right)e_i^\top = P(t) \sum_{i=1}^n e_ie_i^\top = P(t).
\end{equation*}
Similarly,
\begin{equation*}
F^i(t)e_ie_i^\top F^i(t) = P(t) e_i e_i^\top P(t).
\end{equation*}
With these identifications, to show that $(F^1, \ldots, F^n)$ solves \eqref{equation-Fi}, it suffices to check the equations \eqref{ODEsF*}.
For this purpose, let us define $\tau(t) := \tr(P(t))/n > 0$. Omitting the time-dependence, the equation \eqref{ODEsF*} then becomes
\begin{equation*}
- \frac{\tau'}{\tau^2} P e_ie_i^\top P + \frac{1}{\tau}\Dot{P}e_ie_i^\top P + \frac{1}{\tau}Pe_ie_i^\top \Dot{P} - \frac{1}{\tau}P^2 e_ie_i^\top P - \frac{1}{\tau}P e_ie_i^\top P^2 + P e_ie_i^\top P =0,
\end{equation*}
or equivalently, multiplying by $\tau$,
\begin{equation}\label{eq_Pe_i}
- \frac{\tau'}{\tau} P e_ie_i^\top P + \Dot{P}e_ie_i^\top P + P e_ie_i^\top \Dot{P} - P^2e_ie_i^\top P  - Pe_ie_i^\top P^2 + \tau P e_ie_i^\top P =0.
\end{equation}
Hence $(F^1, ..., F^n)$ are solutions of the ODEs \eqref{ODEsF*} if and only if $P$ solves the ODEs \eqref{eq_Pe_i} for $i =1, \ldots, n$, with $\tau(t)=\tr(P(t))/n$.

Now let $v_1,\ldots,v_n$ denote an orthonormal basis of eigenvectors of the symmetric matrix $\lap$, with associated eigenvalues $\lambda_1,\ldots,\lambda_n$. From the definition \eqref{def:P} of $P(t)$ it follows that $v_1,\ldots,v_n$ are eigenvectors of $P(t)$, and the eigenvalue $\rho_j(t)$ of $P(t)$ associated with the eigenvector $v_j$ is given by
\begin{align}\label{eq:Pevector}
\rho_j(t) &= \frac{- f'(T-t) \lambda_j}{1 - f(T-t) \lambda_j} = - \partial_t  \log (1 - f(T-t) \lambda_j). 
\end{align}
Note that because $P(t)$ is symmetric, we have not only $P(t)v_k=\rho_k(t)v_k$ but also $v_k^\top P(t) = \rho_k(t) v_k^\top$ for all $k$.
Now, $P(t)$ satisfies \eqref{eq_Pe_i} for all $i=1,...,n$ if and only if for all $i,j,k$ we have
\begin{equation*}
    v_j^\top \left(-\frac{\tau'}{\tau} P e_ie_i^\top P + \Dot{P}e_ie_i^\top P + P e_ie_i^\top \Dot{P} - P^2e_ie_i^\top P  - Pe_ie_i^\top P^2 + \tau P e_ie_i^\top P\right) v_k =0,
\end{equation*}
which is equivalent to 
\begin{equation*}
    (v_j^\top e_i) (e_i^\top v_k) \left(- \frac{\tau'}{\tau} \rho_j \rho_k + \rho_j'\rho_k + \rho_j\rho_k' - \rho_j^2 \rho_k - \rho_j \rho_k^2 + \tau \rho_j \rho_k \right) = 0,
\end{equation*}
where we omit the time-dependence from $\rho_k=\rho_k(t)$.
Notice that the term in the parenthesis can be written as 
\begin{equation*}
    \left(\rho_k' - \rho_k^2 - \frac{\rho_k}{2}\left( \frac{\tau'}{\tau} - \tau\right) \right) \rho_j + \left(\rho_j' - \rho_j^2 - \frac{\rho_j}{2}\left( \frac{\tau'}{\tau} - \tau\right) \right) \rho_k.
\end{equation*}
To complete the proof, it thus suffices to show that 
\begin{equation*}
\rho_k' - \rho_k^2 - \frac{\rho_k}{2}\left( \frac{\tau'}{\tau} - \tau\right) =0.
\end{equation*}
For $k$ such that $\lap v_k=0$, this holds trivially because $\rho_k \equiv 0$; for all other $k$ we have $\rho_k(t) > 0$ for all $t \in [0,T]$, and by dividing by $\rho_k$, we may show instead that
\begin{equation}\label{eq_rhok_nonzero_solves}
\frac{\rho_k'}{\rho_k} - \rho_k = \frac{1}{2}\left(\frac{\tau'}{\tau} - \tau \right), \ \ \text{ for } k \text{ such that } \rho_k(t) > 0 \ \forall t \in [0,T],
\end{equation}
where $\tau := \frac{1}{n}\tr(P) = \frac{1}{n}\sum_{k=1}^n\rho_k$.

Now, using the form of $\rho_k$ from \eqref{eq:Pevector}, a straightforward computation yields
\begin{align}
\frac{\rho_k'(t)}{\rho_k(t)} - \rho_k(t) = - \frac{f''(T-t)}{f'(T-t)}. \label{pf:rho'/rho-rho}
\end{align}
To simplify the right-hand side of \eqref{eq_rhok_nonzero_solves}, we recall $Q(x):=(\det(I-x\lap))^{1/n}=\prod_{k=1}^n(1-x\lambda_k)^{1/n}$ and compute
\begin{align*}
\tau(t) &= \frac{1}{n}\sum_{k=1}^n\rho_k(t) = - \frac{1}{n}\sum_{k = 1}^n \partial_t  \log (1 - f(T-t)\lambda_k) \\
	&= - \partial_t \log  \prod_{k=1}^n (1 - f(T-t)\lambda_k)^{1/n}, \\
	&= - \frac{\partial}{\partial t} \log  Q(f(T-t)) = \frac{f'(T-t)Q'(f(T-t))}{Q(f(T-t))}.
\end{align*}
Therefore after computing the derivative of $\tau$ and rearranging we obtain
\begin{align}
\frac{\tau'(t)}{\tau(t)} - \tau(t) = - \frac{f''(T-t)}{f'(T-t)} - \frac{f'(T-t)Q''(f(T-t))}{Q'(f(T-t))}. \label{pf:tau'/tau-tau}
\end{align}
Now, since $f$ solves the ODE $f'(t)=cQ'(f(t))$, we have also $f''(t) = c Q''(f(t))f'(t)$, and we compute
\begin{align*}
\frac{f'(T-t)Q''(f(T-t))}{Q'(f(T-t))} &= cQ''(f(T-t)) = \frac{f''(T-t)}{f'(T-t)}.
\end{align*}
Returning to \eqref{pf:tau'/tau-tau}, we deduce
\begin{equation*}
\frac{\tau'(t)}{\tau(t)} - \tau(t) = - 2 \frac{f''(T-t)}{ f'(T-t)}.
\end{equation*}
Recalling \eqref{pf:rho'/rho-rho}, this proves \eqref{eq_rhok_nonzero_solves}, which completes the proof that $(F^1,\ldots,F^n)$ defined in \eqref{pf:def:F*i} solves the desired equations. Thus, the controls $\alpha_i^*(t,\bm x) =  -e_i^\top P(t)\bm x$ form a Nash equilibrium as discussed above.

\subsection{State process dynamics in equilibrium.} \label{se:stateprocess}

The state process in equilibrium is
\begin{equation*}
dX_i(t) = - e_i^TP(t) \bm{X}(t) dt + \sigma dW_i(t),  \qquad i = 1,\ldots,n,
\end{equation*}
which can be written in vector form as 
\begin{equation*}
d\bm{X}(t) = - P(t) \bm{X}(t) dt + \sigma d\bm{W}(t).
\end{equation*}
We can explicitly solve this SDE in terms of $P$, noting that $P(t)$ and $P(s)$ commute for all $t,s \in [0,T]$, to obtain 
\begin{equation*}
\bm{X}(t) = e^{- \int_0^t P(s)ds}\bm{X}(0) +  \sigma  \int_0^t e^{ -\int_s^t P(u)du}d\bm{W}(s).
\end{equation*}
Thus we deduce that in equilibrium the law of $\bm{X}(t)$ is 
\begin{equation*}
    \bm{X}(t) \sim \mathcal{N} \left(e^{- \int_0^t P(s)ds} \bm{X}(0), \sigma^2\int_0^t  e^{-2 \int_s^t P(u)du}ds\right).
\end{equation*}
Using the expression $P(t) = - \partial_t \log(I  - f(T-t)\lap)$, noting that the log of the positive definite matrix is well-defined via power series, we have 
\begin{align*}
\exp\left(-\int_s^tP(u)du\right) &= (I  - f(T-t)\lap)(I  - f(T-s)\lap)^{-1}.
\end{align*}
Hence, $\bm X(t)$ has mean vector
\begin{align*}
\E[\bm X(t)]=(I  - f(T-t)\lap)(I  - f(T)\lap)^{-1}\bm X(0)
\end{align*}
and covariance matrix
\begin{align*}
\Var(\bm X(t))=\sigma^2(I - f(T-t)\lap)^2 \int_0^t (I - f(T-s)\lap)^{-2} ds.
\end{align*}

\subsection{Computing the value of the game.} \label{se:valuecomputation}

We next justify \eqref{def:EQvalue}. Returning to the ansatz \eqref{ansatz} and \eqref{solution-h_i}, and recalling the form of $F^i$ from \eqref{pf:def:F*i}, we find
\begin{align}
\frac{1}{n}\sum_{i=1}^nv_k(t,\bm x) &= \frac{1}{n}\sum_{i=1}^n\left(\frac12 \bm x^\top F^i(t)\bm x + h_i(t) \right) \\
	&= \frac{1}{2\tr(P(t))}\sum_{i=1}^n ( e_i^\top P(t)\bm x)^2 + \frac{\sigma^2}{2n}\sum_{i=1}^n \int_t^T \tr(F^i(s))ds. \label{pf:val1}
\end{align}
To compute the second term, we define $R(x) = \log \frac{1}{n}\sum_{k=1}^n \frac{-\lambda_k }{1 - x\lambda_k}$ for $x \ge 0$, so that
\begin{align*}
\frac{1}{n}\sum_{i=1}^n  \tr(F^i(t)) &= \frac{1}{n}\sum_{i,k=1}^n  e_k^\top F^i(t)e_k  = \frac{1}{\tr(P(t))}\sum_{i,k=1}^n  e_k^\top P(t) e_i e_i^\top P(t)e_k \\
	&= \frac{1}{\tr(P(s))}\sum_{k=1}^n  e_k^\top P(t)^2e_k \\
	&= \frac{\tr(P(t)^2)}{\tr(P(t))} \\
	&= -\sum_{k=1}^n \frac{(f'(T-t))^2\lambda_k^2}{(1-f(T-t)\lambda_k)^2} \Bigg/ \sum_{k=1}^n \frac{f'(T-t)\lambda_k}{1-f(T-t)\lambda_k} \\
	&= -R'(f(T-t))f'(T-t) \\
	&= \partial_t R(f(T-t)) .
\end{align*}
Thus, using $f(0)=0$ and $R(0)=0$, we get
\begin{align*}
\frac{\sigma^2}{2n}\sum_{i=1}^n \int_t^T \tr(F^i(s))ds &= \frac{\sigma^2}{2}[R(f(0)) - R(f(T-t))] = - \frac{\sigma^2}{2}R(f(T-t)).
\end{align*}
Finally, note that
\begin{align*}
R(f(T-t)) = \log \frac{1}{n}\sum_{k=1}^n \frac{-\lambda_k }{1 - f(T-t)\lambda_k} = \log \frac{\tr(P(t))}{nf'(T-t)}.
\end{align*}
Returning to \eqref{pf:val1}, we get
\begin{align*}
\frac{1}{n}\sum_{i=1}^nv_k(t,\bm x) &= \frac{|P(t)\bm x|^2}{2\tr(P(t))} - \frac{\sigma^2}{2}\log \frac{\tr(P(t))}{nf'(T-t)}.
\end{align*}
Plugging in $t=0$ and $\bm x=\bm X^G(0)$, the left-hand side equals $\mathrm{Val}(G)$, and the proof of \eqref{def:EQvalue} is complete.

\subsection{Extension to weighted graphs.} \label{se:extension_general_matrices}
The results of Theorem \ref{thm_equilibrium_n_player} can be extended to weighted graphs, by replacing the matrix $L$ by a more general one. Suppose we have $n$ players $V = \{1, ..., n\}$, and let $L \in \R^{n \times n}$. Consider the following game: Each player $v \in V$, associated with a state process with dynamics \eqref{def:intro:SDE}, wants to minimize the cost 
\begin{equation*}
J_v((\alpha_u)_{u\in V}) = \frac{1}{2} \E \left[\int_0^T |\alpha_v(t, \bm{X}(t))|^2 dt + c \left|e_v^\top L \bm{X}(T) \right|^2 \right].
\end{equation*}
It can be shown that the conclusions of Theorem \ref{thm_equilibrium_n_player} remain true as long as $L$ satisfies the following:
\begin{enumerate}[(i)]
\item Let $\mathrm{A}(L)$ denotes the subgroup of permutations $\varphi$ of $V$  satisfying $L_{\varphi(i) \varphi(j)} = L_{ij}$ for all $i,j \in V$. Then $\mathrm{A}(L)$ acts transitively on $V$. That is, for every $i,j \in V$ there exists $\varphi \in \mathrm{A}(L)$ such that $\varphi(i)=j$.
\item $L$ is symmetric and negative semi-definite. 
\end{enumerate}
The proof is exactly as in Sections \ref{se:ODEanalysis} and \ref{section:n_player_finite}. The first assumption (i) is a natural generalization of the transitivity assumption. In fact, Lemma \ref{le:symmetries} was the only place the structure of a graph Laplacian matrix was really used, and this assumption (i) ensures that Lemma \ref{le:symmetries} remains valid when $\Aut(G)$ is replaced by $\mathrm{A}(L)$.
For examples, note that (i) holds if $L$ is of the form $L_{ij}=w(d(i,j))$, where $w$ is any function and $d$ is the graph distance associated with some transitive graph $G$ on vertex set $V$; then $\mathrm{A}(L) \supset \Aut(G)$.
The second assumption (ii) ensures that the function $Q(x) := (\det(I - xL))^{1/n}$ is smooth and increasing on $\R_+$, so that the ODE $f'(t) = c Q'(f(t))$ is uniquely solvable with $f(0)=0$. Note that $L$ could instead be assumed positive semi-definite, as we can then replace it by $-L$ without affecting the cost function.  Many cases of indefinite matrices $L$ would also work, but these cases require a more careful analysis of the ODE, in light of the singularities of $Q'(x)$.

Theorem \ref{thm:convergence_X^Gn_k} admits a similar extension, as long as one is careful to note that the eigenvalue distributions $\mu_{G_n}$ may now have unbounded support as $n\to\infty$. To ensure that the desired integrals converge, one should assume that $\mu_{G_n} \to \mu$ in a stronger (say, Wasserstein) topology.

\section{A direct but heuristic proof of Theorem \ref{thm_equilibrium_n_player}} \label{se:alternateproof}

In this section we aim to give a more enlightening derivation of the solution given in Theorem \ref{thm_equilibrium_n_player}, compared to the more concise ``guess and check" style of proof presented in Section \ref{section:n_player_finite}. To keep this brief, we will avoid giving full details and instead treat this section as a heuristic argument.
Suppose throughout this section that $G$ is a given finite transitive graph with vertex set $V=\{1,\ldots,n\}$, and omit $G$ from the notations as in $\lap=\lap_G$.

\subsection{More on symmetries.}\label{subs_symmetries}

As a first step, we elaborate on the symmetry discussion of Section \ref{se:symmetries}. Recall the notation $\Aut(G)$ and $R_\varphi$ introduced therein, and note that $\Aut(G)$ acts naturally on $\R^V$ via $(\varphi,\bm x) \mapsto  R_{\varphi} \bm  x$.

Suppose we have a solution $(v_1,\ldots,v_n)$ of the HJB system \eqref{n-player-HJB}. The structure of the game described in Section \ref{se:finitehorizon} is invariant under automorphisms of $G$.  More precisely, this should translate into the following symmetry property for the value functions:
\begin{align}
v_i(t,\bm x) = v_{\varphi(i)}(t, R_{\varphi} \bm x), \qquad i \in V, \ \ \varphi \in \Aut(G), \ \ \bm x \in \R^V. \label{sym-vi}
\end{align}
In particular, if $(v_i(t,\bm x))_{i \in V}$ solves the HJB system \eqref{n-player-HJB}, then a straightforward calculation shows that so does $(v_{\varphi(i)}(t, R_{\varphi} \bm x))_{i \in V}$. Hence, if uniqueness holds for \eqref{n-player-HJB}, then we would deduce \eqref{sym-vi}.
Plugging the quadratic ansatz \eqref{ansatz} into both sides of \eqref{sym-vi} yields
\begin{align*}
\bm x^\top F^i(t) \bm x + h_i(t) = \bm x^\top R_{\varphi}^\top F^{\varphi(i)}(t) R_{\varphi}\bm x + h_{\varphi(i)}(t), \qquad i \in V, \ \ \varphi \in \Aut(G), \ \ \bm x \in \R^V.
\end{align*}
Matching coefficients yields
\begin{align}
F^i(t) = R_{\varphi}^\top F^{\varphi(i)}(t)R_{\varphi}, \qquad h_i(t) = h_{\varphi(i)}(t), \qquad  i \in V, \ \ \varphi \in \Aut(G). \label{pf:keysym1}
\end{align}

This immediately shows that the map $i \mapsto h_i(t)$ is constant on orbits of the action of $\Aut(G)$ on $V$. That is, if $i \in V$ and
\[
\mathrm{Orbit}(i) := \{\varphi(i) : \varphi \in \Aut(G)\},
\]
then $h_k(t) = h_j(t)$ for all $k,j \in \mathrm{Orbit}(i)$. 
Note that the orbits $\{\mathrm{Orbit}(i) : i \in V\}$ form a partition of $V$, and the assumption that $G$ is transitive (Definition \ref{def:transitive}) means precisely that $V$ itself is the only orbit.
Similarly, elaborating on the first identity in \eqref{pf:keysym1}, for any $i,j,k \in V$ we find
\begin{equation}\label{eq:F^i_kj}
e_k^\top F^i(t) e_j = e_k^\top R_{\varphi}^\top F^{\varphi(i)}(t)R_{\varphi} e_j = e_{\varphi(k)}^\top F^{\varphi(i)}(t) e_{\varphi(j)}.
\end{equation}

In the Riccati equation \eqref{equation-Fi}, a key role was played by the matrix $P(t) := \sum_{j = 1}^n F^j(t)e_j e_j^\top$. Under a stronger transitivity assumption on the graph, the symmetry property \eqref{pf:keysym1} is enough to ensure that $P(t)$ is symmetric and commutes with $\lap$ for each $t$:

\begin{definition}
We say that a graph $G$ is \emph{generously transitive} if for each $i,j \in V$ there exists $\varphi \in \Aut(G)$ such that $\varphi(i)=j$ and $\varphi(j)=i$.
\end{definition}

\begin{proposition} \label{P_sym}
If $G$ is generously transitive and \eqref{pf:keysym1} holds, then for each $t \in [0,T]$ the matrix $P(t) := \sum_{j = 1}^n F^j(t)e_j e_j^\top$ satisfies:
\begin{enumerate}[(i)]
\item $P(t) = P(t)^\top$.
\item $P(t)$ and $\lap$ commute.\\
\end{enumerate}
\end{proposition}

\proof{Proof of Proposition \ref{P_sym}.} For brevity, we write $F^k=F^k(t)$ and $P=P(t)$. Note that $Pe_k=F^ke_k$ for each $k=1,\ldots,n$. Let $j,k \in V$, and let $\varphi \in \Aut(G)$ be such that $\varphi(j)=k$ and $\varphi(k)=j$. 
\begin{enumerate}[(i)]
\item We find from \eqref{eq:F^i_kj} that
\begin{equation*}
e_j^\top P e_k = e_j^\top F^k e_k = e_{\varphi(j)}^\top F^{\varphi(k)} e_{\varphi(k)} = e_{k}^\top F^{j} e_j= e_{k}^\top P e_j.
\end{equation*}
\item Using the identities $e_j = R_\varphi e_k$ and $e_k = R_\varphi e_j$, then the fact that $\lap$ commutes with $R_{\varphi}^\top$, and finally that $R_{\varphi}^{\top}F^{k}R_{\varphi} = R_{\varphi}^{\top}F^{\varphi(j)}R_{\varphi} = F^j$ we get
\begin{equation*}
e_j^\top \lap P e_k = e_j^\top \lap F^k e_k = e_k^\top R_{\varphi}^\top \lap F^k R_\varphi e_j = e_k^\top \lap R_{\varphi}^\top F^k R_\varphi e_j = e_k^\top \lap F^j e_j = e_k^\top \lap P e_j.
\end{equation*}
\end{enumerate} 
\endproof

In our ``guess and check" proof of Theorem \ref{thm_equilibrium_n_player} given in Section \ref{section:n_player_finite}, the two properties of Proposition \ref{P_sym} followed automatically from the asserted formula \eqref{def:intro:P*} for $P(t)$. Here we see, on the other hand, that these properties follow from purely algebraic arguments.

\subsection{Heuristic solution of the HJB system.}\label{section_heuristics_thm}

In this section we will explain how to find the expressions of $P$ and $F^i$. We start here from the Riccati equation \eqref{equation-Fi}, and assume from now on that the graph $G$ is generously transitive. The objective is to find $F^1,\ldots, F^n$ such that each $F^i$ is a solution of the matrix differential equation
\begin{equation}\label{Riccati_diff_eq}
   \Dot{F^i}  - P F^i - F^i P^\top +  F^i e_i e_i^\top F^i = 0,
\end{equation}
with terminal condition $F^i(T) = c \lap e_i e_i^\top \lap$, and with $P = \sum_{i=1}^n F^ie_ie_i^\top$. We use a fixed point approach: Treating $P$ as given, this is a system of $n$ decoupled Ricatti differential equations which we can solve. The solutions $F^1, ..., F^n$ give rise to a new $P$, which we then match to the original $P$.

Now fix $P$, which by Proposition \ref{P_sym} we can assume to be symmetric and commuting with $\lap$, and let us solve for $F^1,\ldots,F^n$. From~\cite{reid1963riccati}, we know that if  $(Y^i , \newW^i)$ is solution of the equation
\begin{equation*}
    \begin{bmatrix} \Dot{Y^i} \\ \Dot{\newW^i} \end{bmatrix} = \begin{bmatrix} P & 0 \\  e_ie_i^\top &  -P\end{bmatrix} \begin{bmatrix} Y^i \\ \newW^i\end{bmatrix}.
\end{equation*}
on $[0,T]$ with $\newW^i$ nonsingular on $[0,T]$ then $F^i = Y^i (\newW^i)^{-1}$ is a solution of \eqref{Riccati_diff_eq}. 

Since $P(t)$ is symmetric and commutes with $\lap$, the two matrices are simultaneously diagonalizable. If we further assume that $P(t)$ and $P(s)$ commute for all $s$ and $t$, then the matrices $\lap$ and $\{P(t) : t \in [0,T]\}$ are all simultaneously diagonalizable. We can then choose an orthonormal basis $V=(v_1, \ldots, v_n)$ of eigenvectors of $P(t)$ and $\lap$, such that $\lap$ and $P(t)$ are diagonalizable in this basis for all $t \in [0,T]$. For each $t$, let  $\rho_1(t), \ldots, \rho_n(t)$ denote the eigenvalues of $P(t)$ with corresponding eigenvectors $v_1, \ldots, v_n$.

We can now easily solve the equation for $Y^i$ to get
\begin{equation*}
Y^i(t) =  \Tilde{P}(t) Y^i(T), \qquad \text{ where } \ \ \ \Tilde{P}(t) :=  \exp\left(-\int_t^T P(s)ds\right).
\end{equation*}
We deduce from the expression of $Y^i$ that $\newW^i$ solves the equation $\Dot{\newW}^i(t) = e_ie_i^\top \Tilde{P}(t)Y^i(T) - P(t)\newW^i(t)$, and it follows that
\begin{equation*}
\newW^i(t) = \Tilde{P}^{-1}(t) \left[\newW^i(T) - \int_t^T \Tilde{P}(s)e_ie_i^\top \Tilde{P}(s)Y^i(T)ds \right].
\end{equation*}
If we now choose the boundary values
\begin{equation*}
\newW^i(T) = I, \hspace{1cm} Y^i(T) = c \lap e_i e_i^\top \lap,
\end{equation*}
then $F^i=Y^i(\newW^i)^{-1}$ is symmetric and the terminal condition is satisfied.
Plugging these terminal conditions into $\newW^i$ and $Y^i$ yields
\begin{equation*}
\newW^i(t) = \Tilde{P}^{-1}(t) \left[I - c \int_t^T \frac{\mbox{Tr}(\tilde{P}(s)\lap)}{n} \Tilde{P}(s)e_ie_i^\top\lap ds \right], \quad Y^i(t) = c\Tilde{P}(t)\lap e_ie_i^\top\lap.
\end{equation*}
Here we used the identity $e_i^\top \Tilde{P}(s)\lap e_i = \frac{\tr(\tilde{P}\lap)}{n}$ which holds by Lemma \ref{le:symmetries}(ii), since $\tilde{P}\lap$ commutes with $R_{\varphi}$ for all $\varphi \in \Aut(G)$.

Now that we have our explicit solutions $(Y^i, \newW^i)$, assuming that $\newW^i(t)$ is invertible, we deduce that the Ricatti equation \eqref{Riccati_diff_eq} admits the following solution
\begin{equation}\label{eq:F^i}
F^i(t) =  c\Tilde{P}(t)\lap e_ie_i^\top\lap \left[I - c \newA^i(t) \right]^{-1} \Tilde{P}(t),
\end{equation}
where we define
\begin{align*}
\newA^i(t) := \int_t^T  \frac{\tr(\tilde{P}(s)\lap)}{n}\Tilde{P}(s)e_ie_i^\top\lap ds.
\end{align*}
The objective is now to solve for $P$. To that end, we first note that
\begin{equation}\label{F^i_ei}
F^i(t)e_i =  c \eta_i(t)  \Tilde{P}(t)\lap e_i, \quad \text{where} \quad \eta_i(t) := e_i^\top \lap [I - c  \newA^i(t) ]^{-1} \Tilde{P}(t)e_i . 
\end{equation}

Recalling that both $\lap$ and $\Tilde{P}$ commute with $R_{\varphi}$, a simple computation shows that $R_\varphi \newA^i(t) = \newA^{\varphi(i)}(t)R_\varphi$ for $\varphi \in \Aut(G)$. It follows that $R_\varphi \lap [I - c  \newA^i(t) ]^{-1} \Tilde{P}(t) = \lap [I - c  \newA^{\varphi(i)}(t) ]^{-1} \Tilde{P}(t)R_\varphi$, and using Lemma \ref{le:symmetries}(iii) we deduce that $\eta_1(t)=\eta_2(t)=\cdots=\eta_n(t)$, and we let $\eta(t)$ denote the common value. We then compute
\begin{equation}\label{P_asfct_MPtilde}
P(t) = \sum_{i = 1}^n F^i(t)e_i e_i^\top = c \sum_{i = 1}^n \eta(t) \Tilde{P}(t)\lap e_i e_i^\top = c \eta(t) \tilde{P}(t) \lap.
\end{equation}
Because $P$, $\Tilde{P}$, and $\lap$ are simultaneously diagonalizable, we deduce the following relationship between the eigenvalues:
\begin{equation*}
\rho_k(t) e^{\int_t^T \rho_k(s)ds} = c \eta(t) \lambda_k, \qquad k=1,\ldots,n,
\end{equation*}
where $\lambda_1,\ldots,\lambda_n$ are the eigenvalues of $\lap$, corresponding to eigenvectors $v_1,\ldots,v_n$.
Integrating from $t$ to $T$, taking the logarithm, and finally differentiating leads to
\begin{equation}\label{uk_t}
    \rho_k(t) = \frac{c \lambda_k \eta(t)}{1 + c \lambda_k \int_t^T \eta(s)ds},
\end{equation}
which we can rewrite in matrix form as
\begin{equation}
    P(t) =  c \eta(t) \lap \left(I + c \int_t^T \eta(s)ds \,\lap \right)^{-1}. \label{pf:eq:P}
\end{equation}
This completes our fixed point argument, provided we identify $\eta$. 

Using once more that $e_i^\top\lap \tilde{P}(s)e_i = \frac{\tr(\tilde{P}(s)\lap)}{n}$ for all $s$, a quick computation shows $(\newA^i(t))^2 = \Big(\int_t^T \big(\frac{\tr(\tilde{P}(s)\lap)}{n}\big)^2 ds \Big) \newA^i(t)$. Thus, assuming the validity of the power series $(I-c\newA^i(t))^{-1} = \sum_{k=0}^\infty (c\newA^i(t))^k$, we have
\begin{align}
[I - c \newA^i(t) ]^{-1} &= \sum_{k=0}^{\infty} (c \newA^i(t))^k  = I + c\left[1 - c \int_t^T \Big(\frac{\tr(\tilde{P}(s)\lap)}{n}\Big)^2 ds\right]^{-1} \newA^i(t). \label{pf:Ai-powerseries}
\end{align}
Plugging this back into the definition \eqref{F^i_ei} of $\eta$, we get (for any $i$)
\begin{equation}\label{eq:mu(t)}
\eta(t) = e_i^\top \lap[I-c\newA^i(t)]^{-1}\tilde{P}(t) e_i = \frac{\tr(\tilde{P}(t)\lap)}{n} \, \frac{1}{1 - c \int_t^T \left(\frac{\tr(\tilde{P}(s)\lap)}{n}\right)^2 ds }.
\end{equation}
Using \eqref{P_asfct_MPtilde} and \eqref{uk_t} we get $\tr(\tilde{P}(t)\lap) =  \frac{1}{c \eta(t)} \tr(P(t)) = \sum_{k=1}^n \frac{\lambda_k}{1 + c \lambda_k \int_t^T \eta(s)ds}$, and from \eqref{eq:mu(t)} we then obtain
\begin{equation}\label{g'(t)}
    \eta(t) = \frac{\frac{1}{n} \sum_{k=1}^n \frac{\lambda_k}{1 + c \lambda_k \int_t^T \eta(s)ds}}{1 - c \int_t^T  \Big(\frac{1}{n} \sum_{k=1}^n \frac{\lambda_k}{1 + c \lambda_k \int_s^T \eta(s)ds} \Big)^2 ds}.
\end{equation}
Multiplying both sides by $\frac{1}{n} \sum_{k=1}^n \frac{\lambda_k}{1 + c \lambda_k \int_t^T \eta(s)ds}$, integrating from $t$ to $T$, and exponentiating yield
\begin{equation}\label{geom_mean_g(t)}
    \prod_{k=1}^n \left(1 + c \lambda_k \int_t^T \eta(s)ds \right)^{1/n} = \frac{1}{1 - c \int_t^T  \Big(\frac{1}{n} \sum_{k=1}^n \frac{\lambda_k}{1 + c \lambda_k \int_s^T \eta(s)ds} \Big)^2 ds}.
\end{equation}
Defining $Q(x)=(\det(I - x\lap ))^{1/n}$ and $f(t) = -c\int_{T-t}^T\eta(s)ds$, we find from \eqref{g'(t)} and \eqref{geom_mean_g(t)} that $f$ must satisfy $f(0) = 0$ and $f'(t)=cQ'(f(t))$. Returning to \eqref{pf:eq:P}, and noting that $f'(T-t)=-\eta(t)$, 
we have thus proved that $P$ can be written as
\begin{equation*}
    P(t) = - f'(T-t) \lap \left(I - f(T-t) \lap \right)^{-1},
\end{equation*}
justifying the expression of $P$ in Theorem \ref{thm_equilibrium_n_player}. 

We can also simplify the expression of $F^i$ in \eqref{eq:F^i} to recover the expression we introduced in \eqref{pf:def:F*i} our first the proof of Theorem \ref{thm_equilibrium_n_player}. Use \eqref{pf:Ai-powerseries} to get
\begin{equation*}
F^i(t) = \frac{c}{1 - c \int_t^T\left(\frac{\tr{(\tilde{P}(s)\lap)}}{n}\right)^2 ds } \Tilde{P}(t)\lap e_ie_i^\top\lap \Tilde{P}(t).
\end{equation*}
From (\ref{P_asfct_MPtilde}) and \eqref{eq:mu(t)}, we deduce $1 - c \int_t^T \left(\frac{\tr{(\tilde{P}(s)\lap)}}{n}\right)^2 ds = \frac{1}{c\mu(t)^2}\frac{\tr{(P(t))}}{n}$ and
\begin{equation*}
    F^i(t) = \frac{1}{\tr{(P(t))}/n} P(t) e_ie_i^\top P(t).
\end{equation*}

\section{Correlation decay: Proof of Proposition \ref{pr:corrdecay}}  \label{subsection:cov_bounds}

The purpose of this section is to prove Proposition \ref{pr:corrdecay}, which is essential in the proof of the convergence of the empirical measure given in the next section.

Throughout this section, we fix a finite transitive graph $G = (V,E)$ without isolated vertices and with vertex set $V = \{1, ..., n\}$ for some $n \in \N$. We may without loss of generality assume $G$ to be connected, as otherwise $X^G_v(t)$  and $X_u^G(t)$ are independent for $u$ and $v$ in distinct connected components and the right-hand side of \eqref{def:corrdecaybound} is zero.
Since $G$ is fixed throughout the section we omit it in the notations, e.g., $\lap = \lap_G$, $\bm{X}(t) = \bm{X}^G(t)$, $\delta = \delta(G)$, and $f=f_G$. As before, $(e_1, ..., e_n)$ denotes the standard Euclidean basis in $\R^n$.
We make some use of the adjacency matrix $A=A_G$ in this section, and we repeatedly use the well known fact that $e_u^\top A^\ell e_v$ counts the number of paths of length $\ell$ from $v$ to $u$, for each $\ell \in \N$ and vertices $v,u \in V$.

\proof{Proof of Proposition \ref{pr:corrdecay}.}
Recall from \eqref{intro:variance-lambda} that $X_u(t)$ is Gaussian with variance
\begin{align*}
\frac{\sigma^2}{n}\sum_{k=1}^n \int_0^t  \left(\frac{1 - f(T-t)\lambda_k}{1 - f(T-s)\lambda_k}\right)^2 ds \le \sigma^2 T.
\end{align*}
Indeed, the inequality follows from the fact that $-2 \le \lambda_k \le 0$ and $f$ is increasing.
Since $0 < \gamma < 1$, this proves the claim \eqref{def:corrdecaybound} in the case $u=v$ (noting also that $\delta(G) \ge 1$). We thus focus henceforth on the case of distinct vertices.

From Theorem \ref{thm_equilibrium_n_player}, we know that in equilibrium the state process $\bm{X}$ is normally distributed with covariance matrix
\begin{equation*}
\sigma^2(I - f(T-t)\lap)^2 \int_0^t (I - f(T-s)\lap)^{-2} ds.
\end{equation*}
Our objective is to find a bound for the off-diagonal elements of this matrix. Fix two distinct vertices $u,v\in V$. We have
\begin{equation}\label{eq:cov_1}
\Cov(X_u(t), X_v(t)) = \sigma^2 \int_0^t \sum_{k=1}^n e_u^\top (I - f(T-t)\lap)^2 e_k e_k^\top (I - f(T-s)\lap)^{-2} e_v ds.
\end{equation}
Let us first develop $e_u^\top (I - f(T-t)\lap)^2 e_k$. Using $\lap = \frac{1}{\delta}A - I$, where $A$ is the adjacency matrix of the graph, we find
\begin{equation*}
    \begin{split}
        e_u^\top (I - f(T-t)\lap)^2 e_k = ( 1 + f(T-t))^2 \Bigg(e_u^\top e_k &- 2 \frac{f(T-t)}{\delta( 1 + f(T-t))} e_u^\top Ae_k \\
        & + \left(\frac{f(T-t)}{\delta( 1 + f(T-t))} \right)^2 e_u^\top A^2 e_k\Bigg).
    \end{split}
\end{equation*}
Recall that $d(i,j) = d_G(i,j)$ denotes the distance between two vertices $i$ and $j$ in the graph, i.e., the length of the shortest path from $i$ to $j$ in $G$. Let us write $P(i,m)$ the set of vertices which can be reached in exactly $m$ steps from $i$. By definition of $A$, $e_u^\top A e_k = 1$ if $d(u,k) = 1$ and zero otherwise. Similarly, $e_u^\top A^2 e_k$ is the number of paths of length $2$ from $u$ to $k$, so, in particular, $e_u^\top A^2 e_k=0$ unless $k \in P(u,2)$. Plugging this into \eqref{eq:cov_1}, we get
\begin{align}
\Cov(X_u(t), X_v(t)) &=  \sigma^2 ( 1 + f(T-t))^2 \int_0^t \Bigg\{ e_u^\top (I - f(T-s)\lap)^{-2} e_v \nonumber \\
        &\quad - \sum_{k \in P(u,1)} 2 \frac{f(T-t)}{\delta( 1 +f(T-t))} e_k^T (I - f(T-s)\lap)^{-2} e_v \label{eq:cov_2} \\
        &\quad + \sum_{k \in P(u,2)}  \left(\frac{f(T-t)}{\delta( 1 +f(T-t))} \right)^2 e_u^\top A^2 e_k e_k^T (I - f(T-s)\lap)^{-2} e_v \Bigg\} ds. \nonumber
\end{align}

Next, we estimate the term
\begin{equation*}
    e_k^T (I - f(T-s)\lap)^{-2} e_v =  ( 1 +f(T-s))^{-2} e_k^T \left(I - \frac{f(T-s)}{\delta( 1 +f(T-s))}A\right)^{-2} e_v.
\end{equation*}
To simplify the notations, we define the function
\begin{equation*}
    \gamma(s) := \frac{f(T-s)}{ 1 +f(T-s)}.
\end{equation*}
From  Proposition \ref{pr:f-convergence}, we know that $0 \le f(t)\le ct$ for all $t \in [0,T]$, and thus $0 \le \gamma(s) < 1$ for all $s \in [0,T]$. Moreover, the spectral radius of the adjacency matrix $A$ is always bounded by the degree $\delta$. We can thus use the series $\frac{1}{(1-x)^2} = \sum_{\ell=0}^{\infty} (\ell+1)x^\ell$ for $|x| < 1$ to get
\begin{equation}
    \begin{split}
        e_k^\top (I - f(T-s)\lap)^{-2} e_v &= ( 1 +f(T-s))^{-2} \sum_{\ell = 0}^{\infty} (\ell+1) \left( \frac{\gamma(s)}{\delta} \right)^\ell e_k^\top A^\ell e_v \\
        &= ( 1 +f(T-s))^{-2} \sum_{\ell = d(k,v)}^{\infty} (\ell+1) \left( \frac{\gamma(s)}{\delta } \right)^\ell e_k^\top A^\ell e_v,
    \end{split} \label{pf:exp-2bound}
\end{equation}
where in the last line we used the fact that $e_k^T A^\ell e_v = 0$ if $\ell < d(k,v)$, since the latter implies there are no paths of length $\ell$ from $k$ to $v$.
Next, note the elementary bound $0 \le e_k^T A^\ell e_v \le \delta^{\ell - 1}$, since each vertex has exactly $\delta$ neighbors. Indeed, for each of the first $\ell-1$ steps of a path of length $\ell$ from $k$ to $v$, there are at most $\delta$ choices of vertex, and for the last step there is at most one choice which terminates the path at $v$. Hence, using the identity
\begin{equation*}
\sum_{\ell=k}^{\infty} (\ell+1)x^\ell  = \frac{d}{d x} \left(\sum_{\ell=k}^{\infty} x^{\ell+1} \right) = \frac{d}{d x} \left(\frac{x^{k+1}}{1 - x} \right) = \frac{x^k}{(1-x)^2} (1 + k(1-x))
\end{equation*}
for $|x| < 1$, we deduce from \eqref{pf:exp-2bound} that
\begin{equation*}
   0 \le  e_k^\top (I - f(T-s)\lap)^{-2} e_v \le \frac{1}{\delta( 1 +f(T-s))^{2}}  \frac{\gamma(s)^{d(k,v)} (1 + d(k,v)(1 - \gamma(s))}{(1 - \gamma(s))^2}.
\end{equation*}
Now notice that in \eqref{eq:cov_2} the first term and all the terms in the last sum on the right hand side are positive, whereas all the terms in the second line sum are negative. Therefore we have the following upper bound
\begin{align}
\Cov(X_u(t), X_v(t)) \le \, &\sigma^2 \int_0^t \frac{(1 +f(T-t))^2}{\delta( 1 +f(T-s))^{2}}  \Bigg\{ \frac{\gamma(s)^{d(u,v)} (1 + d(u,v)(1 - \gamma(s)))}{(1 - \gamma(s))^2} \label{eq:cov_upper_bound} \\
     &+ \sum_{k \in P(u,2)} e_u^\top A^2 e_k \left(\frac{\gamma(t)}{\delta} \right)^2   \frac{\gamma(s)^{d(k,v)} (1 + d(k,v)(1 - \gamma(s)))}{(1 - \gamma(s))^2} \Bigg\} ds. \nonumber
\end{align}
The function $d \mapsto \gamma(t)^d(1 + d(1 - \gamma(t)))$ is non-increasing on $[0, \infty)$ since $\gamma(t) \in [0,1)$. Moreover, for all $k$ such that $d(k,u) \le 2$, we have $|d(k,v) - d(u,v)| \le 2$, and in particular $d(k,v) \ge d(u,v) - 2$. We consider two cases separately:
\begin{enumerate}
\item First, suppose $d(u,v) \ge 2$. Since $d(k,v) \ge d(u,v)-2$ by the above argument, monotonicity of $d \mapsto \gamma(t)^d(1 + d(1 - \gamma(t)))$ lets us estimate
\begin{align*}
\sum_{k \in P(u,2)} & e_u^\top A^2 e_k \left(\frac{\gamma(t)}{\delta} \right)^2   \frac{\gamma(s)^{d(k,v)} (1 + d(k,v)(1 - \gamma(s)))}{(1 - \gamma(s))^2} \\
	&\le  \left(\frac{\gamma(t)}{\delta} \right)^2   \frac{\gamma(s)^{d(u,v)-2} (1 + (d(u,v)-2)(1 - \gamma(s)))}{(1 - \gamma(s))^2}\sum_{k \in P(u,2)} e_u^\top A^2 e_k \\
	&\le \gamma(t)^2   \frac{\gamma(s)^{d(u,v)-2} (1 + (d(u,v)-2)(1 - \gamma(s)))}{(1 - \gamma(s))^2}.
\end{align*}
Indeed, the last step uses the fact that $\sum_{k \in P(u,2)} e_u^\top A^2 e_k$ is precisely the number of paths of length two originating from vertex $u$, which is clearly bounded from above by $\delta^2$.
Moreover, since $f$ is increasing by Proposition \ref{pr:fQproperties}, it is straightforward to check that $\gamma$ is decreasing. Therefore $0 \le \gamma(t) \le \gamma(s) \le 1$ for all $t \ge s \ge 0$, and the above quantity is further bounded from above by
\begin{align*}
\frac{\gamma(s)^{d(u,v)} (1 + d(u,v)(1 - \gamma(s)))}{(1 - \gamma(s))^2}.
\end{align*}
Plugging this back in \eqref{eq:cov_upper_bound}, and using the inequality $\frac{1 +f(T-t)}{1 +f(T-s)} \le 1$ which again follows from the fact that $f$ is increasing and nonnegative, we get
\begin{equation}\label{eq:Cov_ub}
\Cov(X_u(t), X_v(t)) \le \frac{2\sigma^2}{\delta}  \int_0^t \frac{\gamma(s)^{d(u,v)}(1 + d(u,v)(1 - \gamma(s)))}{(1 - \gamma(s))^2}  ds.
\end{equation}
\item Suppose next that $d(u,v)=1$. We then use the bounds $\gamma(t)^d(1 + d(1 - \gamma(t))) \le 1$ and $\sum_{k \in P(u,2)} e_u^\top A^2 e_k \le \delta^2$ to estimate
\begin{align*}
\sum_{k \in P(u,2)} & e_u^\top A^2 e_k \left(\frac{\gamma(t)}{\delta} \right)^2   \frac{\gamma(s)^{d(k,v)} (1 + d(k,v)(1 - \gamma(s)))}{(1 - \gamma(s))^2} \\
	&\le \frac{\gamma(t)^2}{(1 - \gamma(s))^2} \le \frac{\gamma(s)^2}{(1 - \gamma(s))^2} \le \frac{\gamma(s)^{d(u,v)+1}(1+d(u,v)(1-\gamma(s))}{(1 - \gamma(s))^2}.
\end{align*}
Plugging this into \eqref{eq:cov_upper_bound} shows that the same bound \eqref{eq:Cov_ub} is valid for $d(u,v)=1$.
\end{enumerate}

Now, recall that $\gamma(\cdot)$ is decreasing. Since $f(T)\le cT$ by Proposition \ref{pr:f-convergence}, we have $0 < \gamma(0) = \frac{cf(T)}{1+cf(T)} \le \frac{cT}{1 + cT} < 1$. Set $\gamma := \frac{cT}{1 + cT}$. The function $y \mapsto \frac{y^d(1+d(1-y))}{(1 - y)^2}$ is easily seen to be increasing on $[0,1)$ for any $d \ge 0$. Thus, from \eqref{eq:Cov_ub} we finally deduce the desired upper bound
\begin{equation*}
\Cov(X_u(t), X_v(t)) \le 2\sigma^2 t   \frac{ \gamma^{d(u,v)}(1 + d(u,v)(1-\gamma))}{\delta(1 - \gamma)^2}.
\end{equation*}

Using the same arguments for the second term in \eqref{eq:cov_2}, which is the only negative term, we obtain a similar lower bound for $\Cov(X_u(t), X_v(t))$, which concludes the proof. 
\endproof

\begin{remark} \label{re:decayofcontrol}
The arguments given here could likely be adapted to estimate the dependence of the equilibrium control of one player on a distant player's state, which would provide an interesting alternative notion of ``correlation decay". To be precise, recall the equilibrium control $\alpha=\alpha^G$ from Theorem \ref{thm_equilibrium_n_player}. For two vertices $(i,j)$, we have  $\partial_{x_j}\alpha_i(t,\bm x) = f'(T-t) e_i^\top L(I-f(T-t)L)^{-1} e_j$, and we suspect that similar arguments to those given above could show that this matrix entry decays exponentially with the graph distance $d(i,j)$. We do not pursue this, as it is not directly suited to our application to empirical measure convergence of the next section.
\end{remark}

\section{Asymptotic regimes}\label{se:asymptoticregimes}

In this section, we provide the derivations of the large-$n$ asymptotics of the in-equilibrium processes. We will first prove Theorem \ref{thm:convergence_X^Gn_k}, then Theorem \ref{th:approxeq}, and lastly we will focus on the examples we discussed in Section \ref{se:examples}, and in particular prove Corollary \ref{co:dense} and Proposition \ref{prop:cyclegraph}.

\subsection{Large-scale asymptotics on transitive graphs: Proof of Theorem \ref{thm:convergence_X^Gn_k}.}

Part (1) of Theorem \ref{thm:convergence_X^Gn_k} is a consequence of Proposition \ref{pr:f-convergence}, so we focus on parts (2--4).
Let $\{G_n\}$ be a sequence of finite transitive graphs, and let $\{\mu_{G_n}\}$ be the corresponding sequence of empirical eigenvalue distributions defined by \eqref{def:spectralmeasure}. We assume that $\{\mu_{G_n}\}$ converges weakly to a probability measure $\mu$. Recall that the initial states are $\bm{X}^{G_n}(0) = \bm{0}$, and recall from \eqref{intro:variance-lambda} that each $X^{G_n}_i(t)$ is Gaussian with mean zero and variance
\begin{align*}
\varfunc_{G_n}(t) &:=  \sigma^2 \int_0^t \int_{[-2,0]} \left(\frac{1 - f_{G_n}(T-t)\lambda}{1 - f_{G_n}(T-s)\lambda}\right)^2 \mu_{G_n}(d\lambda)ds.
\end{align*}

\subsubsection{Convergence of $X^{G_n}_{k_n}(t)$: Proof of (2).}

By Proposition \ref{pr:f-convergence}, we know $f_{G_n}=f_{\mu_{G_n}}$ converges uniformly to the function $f_\mu$ given by \eqref{def:limitODE}. Defining $\varfunc_\mu(t)$ as in \eqref{def:gaussianlimit}, it follows from this uniform convergence and the weak convergence of $\mu_{G_n}$ to $\mu$ that $\varfunc_{G_n}(t) \to \varfunc_\mu(t)$. Therefore $X^{G_n}_{k_n}(t) \sim\varfunc_{G_n}(t)$ converges weakly to $\mathcal{N}(0,\varfunc_\mu(t))$ as $n\to\infty$.

\subsubsection{Convergence of the empirical measure: Proof of (3).} \label{se:convempmeas}

We next show that the (random) empirical measure
\begin{align*}
m^{G_n}(t) := \frac{1}{|G_n|}\sum_{v \in G_n} \delta_{X^{G_n}_v(t)}
\end{align*}
converges to the Gaussian measure $\mathcal{N}(0,\varfunc_\mu(t))$, for each $t \in [0,T]$. In fact, it suffices to show that $m^{G_n}(t)$ concentrates around its mean, in the following sense: For any bounded 1-Lipschitz function $h$, it holds that
\begin{equation}
\lim_{n\to\infty}\E \left[\left|\int h \, dm^{G_n}(t) - \E \int h \, dm^{G_n}(t)\right|^2\right] = 0. \label{pf:empconcentration}
\end{equation}
Indeed, once \eqref{pf:empconcentration} is established, it follows from the transitivity of $G_n$ that $\E \int h \, dm^{G_n}(t) = \E[h(X^{G_n}_{k_n}(t))]$, where $k_n \in G_n$ is arbitrary. Since the law of $X^{G_n}_{k_n}(t)$ converges weakly to $m(t):=\mathcal{N}(0,\varfunc_\mu(t))$, we deduce that $\int h \, dm^{G_n}(t) \to \int h \, dm(t)$ in probability, and the claim follows.

Before proving \eqref{pf:empconcentration}, we digress to state a lemma pertaining to the degrees. Recall that each vertex in the transitive graph $G_n$ has the same degree, denoted $\delta(G_n)$.

\begin{lemma}\label{lemma_degree}
We have $\mu_{G_n} \to \delta_{-1}$ if and only if $\delta(G_n) \to \infty$.
If $\mu_{G_n} \to \mu \neq \delta_{-1}$, then $\sup_n \delta(G_n) < \infty$.\\
\end{lemma}
\proof{Proof of Lemma \ref{lemma_degree}.}
Recall that $\lap_{G_n}=\frac{1}{\delta(G_n)}A_{G_n}- I$, where $A_{G_n}$ is the adjacency matrix of the graph $G_n$. Then
\begin{align*}
\Var(\mu_{G_n}) &= \int_{[-2,0]} (1+\lambda)^2 \,\mu_{G_n}(dx) = \frac{1}{n}\tr\Big(\frac{1}{\delta(G_n)^2}A_{G_n}^2\Big) = \frac{1}{n\delta(G_n)^2}\sum_{i=1}^n (A_{G_n}^2)_{ii}.
\end{align*}
Since $(A_{G_n}^2)_{ii}=\delta(G_n)$ is exactly the number of paths of length $2$ starting and ending at vertex $i$, we get $\Var(\mu_{G_n}) = 1/\delta(G_n)$. 
Thus, if $\mu_{G_n} \to \mu$ weakly for some probability measure $\mu$ on $[-2,0]$, we have
\begin{align}
\Var(\mu)=\int_{[-2,0]} (1+\lambda)^2 \,\mu(dx) &= \lim_{n\to\infty}\frac{1}{\delta(G_n)}. \label{second moment}
\end{align}
The second claim follows immediately.
It is straightforward to check that $\mu_{G_n} \to \delta_{-1}$ if and only if $\Var(\mu_{G_n}) \to 0$, and the first claim follows.
\endproof

We now turn toward the proof of \eqref{pf:empconcentration}, for a fixed bounded 1-Lipschitz function $h$. We achieve this by applying the Gaussian Poincar\'e inequality and then using the covariance estimate of \ref{pr:corrdecay}. To this end, fix $t \in [0,T]$ and $n$, and suppose for simplicity that $G_n$ has $n$ vertices. Let $\Sigma^n$ denote the $n \times n$ covariance matrix of $\bm X^{G_n}(t)$, and let $M$ denote its symmetric square root. Then $\bm X^{G_n}(t) \stackrel{d}{=} M \bm Z$ for a standard Gaussian $\bm Z$ in $\R^n$. Define $F : \R^n \to \R$ by
\begin{align*}
F(\bm x) := \frac{1}{n}\sum_{i=1}^n h(e_i^\top M \bm x).
\end{align*}
Then $F(\bm Z) \stackrel{d}{=} \int h\,dm^{G_n}(t)$. 
Noting that $\partial_i F(x) = \frac{1}{n} \sum_{j=1}^n h'(e_j^\top M \bm x) M_{ji}$, we get
\begin{equation*}
    \begin{split}
        |\nabla F(\bm x)|^2 &= \sum_{i=1}^n \frac{1}{n^2} \Big(\sum_{j=1}^n h'(e_j^\top  M \bm x) M_{ji}\Big)^2 \le \frac{1}{n^2} \sum_{j,k = 1}^n h'(e_j^\top M \bm x) h'(e_k^\top M \bm x) \Sigma_{jk} \\
        & \le \frac{1}{n^2} \sum_{j,k = 1}^n |\Sigma_{jk}^n|,
    \end{split}
\end{equation*}
where in the last inequality we used the fact that $h$ is $1$-Lipschitz. Now, applying the Gaussian Poincar\'e inequality (see \cite[Theorem 3.20]{boucheron2013concentration}), we find
\begin{align}
\E \left[\left|\int h \, dm^{G_n}(t) - \E \int h \, dm^{G_n}(t)\right|^2\right] &= \Var(F(\bm Z)) \le \E[|\nabla F(\bm Z)|^2] \le \frac{1}{n^2} \sum_{j,k = 1}^n |\Sigma_{jk}^n|. \label{pf:avgcovar}
\end{align}
It remains to show that this converges to zero as $n\to\infty$.

Let $\epsilon > 0$. By Proposition \ref{pr:corrdecay}, $|\Sigma_{jk}^n|$ converges to $0$ as $d_{G_n}(j,k) \rightarrow \infty$, where $d_{G_n}$ denotes the graph distance in $G_n$.
Choose $m \in \N$ large enough so that $|\Sigma_{jk}^n| \le \epsilon$ for all $n \in \N$ and $j,k \in G_n$ with $d_{G_n}(j,k) > m$.
For $k \in G_n$ let $B_n(j,m)$ denote the set of vertices in $G_n$ of distance at most $m$ from $j$. 
Because $G_n$ is transitive, the cardinality $|B_n(j,m)|$ does not depend on $j \in G_n$, and we denote by $|B_n(m)|$ this common value.
Then, we use the bound on $|\Sigma^n_{jk}|$ from Proposition \ref{pr:corrdecay} to get
\begin{equation*}
    \begin{split}
        \frac{1}{n^2} \sum_{j,k = 1}^n |\Sigma^n_{jk}| &\le \epsilon + \frac{1}{n^2} \sum_{j=1}^n \sum_{k \in B_n(j,m)} |\Sigma^n_{jk}| \le \epsilon +  \frac{|B_n(m)|}{n\delta(G_n)} \frac{2\sigma^2T}{(1-\gamma)^2},
    \end{split}
\end{equation*}
where we used the fact that $\gamma^d(1+d(1-\gamma))$ is decreasing in $d \ge 0$ since $0 \le \gamma <1$ and is thus bounded by $1$.
We now distinguish two cases. If $\delta(G_n) \to \infty$, then we use $|B_n(m)|/n \le 1$ to send $n\to\infty$ and then $\epsilon \to 0$ to get that \eqref{pf:avgcovar} converges to zero. On the other hand, suppose $\delta(G_n)$ does not converge to infinity. Then necessarily $\sup_n\delta(G_n) < \infty$ by Lemma \ref{lemma_degree}, since $\mu_{G_n}$ converges weakly by assumption. 
Using the obvious bound $|B_n(m)| \le \delta(G_n)^m$, we can again send $n\to\infty$ and then $\epsilon \to 0$ to get that \eqref{pf:avgcovar} converges to zero. This completes the proof of part (3) of Theorem \ref{thm:convergence_X^Gn_k}.

\subsubsection{Convergence of the value: Proof of (4).}

Recall the identity for the value of the game from \eqref{def:valueidentity}. Since $f_{G_n}(T) \to f_\mu(T)$ by Proposition \ref{pr:f-convergence} and $\mu_{G_n} \to \mu$ weakly, it follows that
\begin{equation*}
 \mbox{Val}(G_n) = - \frac{\sigma^2}{2} \log  \int_{[-2,0]} \frac{- \lambda}{1 - f_{G_n}(T)\lambda} \mu_{G_n}(d\lambda) \to - \frac{\sigma^2}{2} \log  \int_{[-2,0]} \frac{- \lambda}{1 - f_\mu(T)\lambda} \mu(d\lambda).
\end{equation*}
This gives (4) and completes the proof of Theorem \ref{thm:convergence_X^Gn_k}.

\subsection{Approximate equilibria: Proof of Theorem \ref{th:approxeq}.}

We begin toward proving Theorem \ref{th:approxeq} by first studying the control $\alpha^{\mathrm{MF}}$ introduced therein. The following lemma shows that it arises essentially as the equilibrium of a mean field game, or equivalently as the optimal control for a certain control problem:

\begin{lemma}\label{lemma_alpha*_optimal}
Define $\alpha^{\mathrm{MF}} : [0,T] \times \R \to \R$ 
\begin{equation}
\alpha^{\mathrm{MF}}(t,x) := - \frac{cx}{1+c(T-t)}. \label{pf:def:optMFcontrol}
\end{equation}
Let $(\Omega',\F',\FF',\PP')$ be any filtered probability space supporting an $\FF'$-Brownian motion $W$ and a $\FF'$-progressively measurable real-valued process $(\beta(t))_{t \in [0,T]}$  satisfying $\E\int_0^T\beta(t)^2dt < \infty$.
Let $X$ be the unique strong solution of the SDE
\begin{align*}
dX(t) &= \alpha^{\mathrm{MF}}(t,X(t))dt + \sigma dW(t), \quad X(0)=0,
\end{align*}
and define $(Y(t))_{t \in [0,T]}$ by
\begin{align*}
dY(t) &= \beta(t)dt + \sigma dW(t), \quad Y(0)=0.
\end{align*}
Then 
\begin{align*}
\frac12\E\left[\int_0^T|\alpha^{\mathrm{MF}}(t,X(t))|^2dt + c|X(T)|^2\right] \le \frac12\E\left[\int_0^T|\beta(t)|^2dt + c|Y(T)|^2\right]
\end{align*}
\end{lemma}
\proof{Proof of Lemma \ref{lemma_alpha*_optimal}.}
We study the HJB equation corresponding to this control problem, which is
\begin{equation*}
\partial_t v(t,x) + \inf_{a \in \R} \left(a  \partial_x v(t,x) + \frac{1}{2} a^2 \right) + \frac{1}{2} \sigma^2 \partial_{xx} v(t,x) = 0,
\end{equation*}
or equivalently
\begin{equation*}
\partial_t v(t,x) - \frac12 |\partial_x v(t,x)|^2 + \frac{1}{2} \sigma^2 \partial_{xx} v(t,x) = 0,
\end{equation*}
with terminal condition $v(T,x) = cx^2/2$.
The ansatz $v(t,x)=a(t)x^2 + b(t)$ yields a classical solution, where $a$ and $b$ are functions satisfying
\begin{align*}
a'(t) - 2 a(t)^2 = 0, \qquad b'(t) + \sigma^2 a(t) = 0,
\end{align*}
with terminal conditions $a(T) = c/2$ and $b(T)=0$. We deduce 
\begin{align*}
a(t) = \frac{1}{2}\frac{c}{1+c(T-t)}, \quad        b(t) = \frac{\sigma^2}{2}\log(1+c(T-t)).
\end{align*}
Therefore the optimal Markovian control is $- \partial_x v(t,x) = - 2 a(t) x = \alpha^{\mathrm{MF}}(t,x)$.
This completes the proof, by a standard verification argument (see \cite[Theorem 3.5.2]{pham2009continuous}).
\endproof

Now let $G$ be a fixed finite graph with vertex set $V=\{1,\ldots,n\}$. Let us again omit the $G$ superscripts from the notation, with $\A=\A_G$ and $J=J^G$ denoting the control set and value function from Section \ref{se:finitehorizon}.
Let $\bm\alpha= (\alpha^{\mathrm{MF}}_i)_{i=1}^n$, and for $\beta \in \A$ and $i \in V$ let $(\beta,\bm\alpha^{-i}) := (\alpha^{\mathrm{MF}}_1,\ldots,\alpha^{\mathrm{MF}}_{i-1},\beta,\alpha^{\mathrm{MF}}_{i+1},\ldots,\alpha^{\mathrm{MF}}_n)$.
To be clear about the notation, we write $\alpha^{\mathrm{MF}}$ without a subscript to denote the control in \eqref{pf:def:optMFcontrol}, whereas $\alpha^{\mathrm{MF}}_i(t,\bm x)=\alpha^{\mathrm{MF}}(t,x_i)$ for $i \in V$ and $\bm x \in \R^n$.
If $v \in V$ has $\deg_G(v) = 0$, then (recalling the definition \eqref{def:cost-isolated} of the cost function for isolated vertices) Lemma \ref{lemma_alpha*_optimal} ensures that 
\begin{align*}
J_v(\bm\alpha) = \inf_{\beta \in \A}J_v(\beta,\bm\alpha^{-v}),
\end{align*}
so we may take $\epsilon_v=0$ as claimed in Theorem \ref{th:approxeq}.
Thus, we assume henceforth that $v \in G$ is a fixed non-isolated vertex, so that $\deg_G(v) \ge 1$.

Now define 
\begin{equation*}
\beta := \arg\!\min_{\beta' \in \A} J_v(\beta', \bm\alpha^{-v}).
\end{equation*}
(Or take a $\delta$-optimizer in case no optimizer exists, and  send $\delta \to 0$ at the end of the proof.)
We aim to prove that
\begin{align}
J_v(\beta,\bm\alpha^{-v}) &\ge J_v(\bm\alpha) -  \frac{\sigma^2 cT}{1+cT} \sqrt{ \frac{cT(2+cT)}{\deg_G(v)}}. \label{pf:approxEQineq}
\end{align}
Define the state processes $\bm X=(X_1,\ldots,X_n)$ and $\bm Y=(Y_1,\ldots,Y_n)$ as the unique strong solutions of the SDEs
\begin{align}
dX_i(t) &= \alpha^{\mathrm{MF}}(t, X_i(t))dt  + \sigma dW_i(t), \quad X_i(0)=0, \ \ i \in V, \label{pf:Xsde} \\
dY_i(t) &= \alpha^{\mathrm{MF}}(t, Y_i(t))dt  + \sigma dW_i(t), \quad Y_i(0)=0, \ \ i \in V \setminus \{v\}, \nonumber \\
dY_v(t) &= \beta(t, \bm Y(t))dt + \sigma dW_v(t), \quad Y_v(0)=0. \nonumber
\end{align}
Note that $Y_i\equiv X_i$ for $i \neq v$.
The values for the player $v$, under $\alpha^{\mathrm{MF}}$ and the deviation $\beta$, are then, respectively,
\begin{align*}
J_v(\bm\alpha) &= \frac{1}{2} \E \left[ \int_0^T |\alpha^{\mathrm{MF}}(t,X_v(t))|^2dt + c \left|\frac{1}{\deg_G(v)}\sum_{u \sim v} X_u(T) - X_v(T)\right|^2\right], \\
J_v(\beta, \bm\alpha^{-v}) &= \frac{1}{2} \E \left[ \int_0^T |\beta(t,\bm{Y}(t))|^2dt + c \left|\frac{1}{\deg_G(v)}\sum_{u \sim v} X_u(T) - Y_v(T)\right|^2\right],
\end{align*}
where we recall that $u \sim v$ means that $u$ is adjacent to $v$.
We prove \eqref{pf:approxEQineq} in three steps:
\begin{enumerate}
\item We show that 
\begin{align*}
J_v(\beta, \bm\alpha^{-v}) \ge \ &\frac{1}{2} \E \left[ \int_0^T |\beta(t, \bm{Y}(t))|^2dt + c |Y_v(T)|^2\right] + \frac{c\sigma^2T}{2\,\deg_G(v)(1+cT)} \\
	&-  c \sqrt{\frac{\sigma^2 T^2}{\deg_G(v)(1+ cT)} \E\int_0^T |\beta(t,\bm{Y}(t))|^2 dt }.
\end{align*}
\item We then estimate
\[
\E \int_0^T |\beta(t,\bm{Y}(t))|^2 dt \le c\sigma^2T\frac{2+cT}{1+cT}.
\]
\item Finally, we show that
\begin{align*}
\frac{1}{2} \E \left[ \int_0^T |\beta(t,\bm{Y}(t))|^2dt + c |Y_v(T)|^2\right] \ge J_v(\bm{\alpha}^{MF}) - \frac{c\sigma^2T}{2\,\deg_G(v)(1+cT)}.
\end{align*}
which will conclude the proof.
\end{enumerate}

\noindent\textbf{Step 1.} We start with 
\begin{align*}
& J_v(\beta, \bm\alpha^{-v}) - \frac{1}{2} \E \left[\int_0^T |\beta(t, \bm{Y}(t))|^2dt + c |Y_v(T)|^2\right] \\
        &\quad = \frac{c}{2} \E \left[ \left|\frac{1}{\deg_G(v)}\sum_{u \sim v} X_u(T) - Y_v(T)\right|^2  -  |Y_v(T)|^2\right] \\
        &\quad = \frac{c}{2} \E \left[ \left(\frac{1}{\deg_G(v)}\sum_{u \sim v} X_u(T)\right)^2\right] - c \E \left[\frac{1}{\deg_G(v)}\sum_{u \sim v} X_u(T) Y_v(T) \right].
\end{align*}
From the form of the SDE \eqref{pf:Xsde} and the definition of $\alpha^{\mathrm{MF}}$, we find that $(X_u(T))_{u \in V}$ are i.i.d.\ Gaussians with mean zero and variance $\sigma^2 T/(1+cT)$. We deduce
\begin{equation}
\E \left[ \left(\frac{1}{\deg_G(v)}\sum_{u \sim v} X_u(T)\right)^2\right] = \frac{\sigma^2T}{\deg_G(v)(1+cT)}. \label{pf:avgvar1}
\end{equation}
For the second term, we note
\begin{equation*}
\E\left[\frac{1}{\deg_G(v)}\sum_{u \sim v} X_u(T) Y_v(T)\right] = \E \left[\frac{1}{\deg_G(v)}\sum_{u \sim v} X_u(T) \int_0^T \beta(t,\bm{Y}(t)) dt\right],
\end{equation*}
and we use Cauchy-Schwarz to bound this term in absolute value by
\begin{equation*}
\sqrt{\frac{\sigma^2T^2}{\deg_G(v)(1+cT)} \E \int_0^T |\beta(t,\bm{Y}(t))|^2 dt }.
\end{equation*}
Plugging these two terms back in our first inequality we obtain claim (1) above.

{\ }

\noindent\textbf{Step 2.} To find a bound for $\E\int_0^T|\beta(t,\bm{Y}(t))|^2dt$, we use the definition of $\beta$ as the minimizer of $J_v(\beta,\bm\alpha^{-v})$. In particular, $J_v(\beta,\bm\alpha^{-v}) \le J_v(0,\bm\alpha^{-v})$. 
Expand this inequality, discarding the non-negative terminal cost on the left-hand side, and noting that the state process of player $v$ when adopting the zero control is precisely $(\sigma W_v(t))_{t \in [0,T]}$, to get
\begin{align*}
\E \int_0^T |\beta(t,\bm{Y}(t))|^2dt &\le c \E \left[  \Big|\frac{1}{\deg_G(v)}\sum_{u \sim v} X_u(T) - \sigma W_v(T)\Big|^2\right].
\end{align*}
Using \eqref{pf:avgvar1} and the independence of the $X_u(T)$ and $W_v(T)$, the right-hand side equals
\begin{align*}
 \frac{c\sigma^2T}{\deg_G(v)(1+cT)} + c\sigma^2 T \le c\sigma^2T\frac{2+cT}{1+cT}.
\end{align*}

\noindent\textbf{Step 3.}
We use Lemma \ref{lemma_alpha*_optimal} to deduce that
\begin{align*}
\frac{1}{2} \E \left[ \int_0^T |\beta(t,\bm{Y}(t))|^2dt + c |Y_v(T)|^2\right] &\ge \frac{1}{2} \E \left[ \int_0^T |\alpha^{\mathrm{MF}}(t,X_v(t))|^2dt + c |X_v(T)|^2\right].
\end{align*}
The right-hand side can be written as
\begin{align*}
J_v(\bm\alpha) + \frac{c}{2}\E\left[ |X_v(T)|^2 - \Big|\frac{1}{\deg_G(v)}\sum_{u \sim v} X_u(T) -  X_v(T) \Big|^2\right].
\end{align*}
Since $(X_u(T))_{u \in V}$ are i.i.d.\ with mean zero and variance $\sigma^2 T/(1+cT)$, as in \eqref{pf:avgvar1} we get 
\begin{align*}
\frac{c}{2}\E\left[ |X_v(T)|^2 - \Big|\frac{1}{\deg_G(v)}\sum_{u \sim v} X_u(T) -  X_v(T) \Big|^2\right] = - \frac{c \sigma^2 T}{2 \deg_G(v) (1+cT)}.
\end{align*}

\subsection{Examples.} \label{se:examples-proofs}

We next specialize the results to the examples discussed in Section \ref{se:examples}.
In particular, we prove Corollary \ref{co:dense} and Proposition \ref{prop:cyclegraph}. 

\subsubsection{Dense case: Proof of Corollary \ref{co:dense}.}\label{section_dense_case}

As in Corollary \ref{co:dense}, let $\{G_n\}$ be a sequence of transitive graphs such that each vertex of $G_n$ has common degree $\delta(G_n) \ge 1$.
The claim that $\mu_{G_n} \to \delta_{-1}$ if and only if $\delta(G_n) \to \infty$ holds as a consequence of Lemma \ref{lemma_degree}.
In this case, we apply Theorem \ref{thm:convergence_X^Gn_k} with $\mu=\delta_{-1}$. The function $Q_\mu$ therein is then $Q_{\delta_{-1}}(x) = 1+x$, and the function $f_{\delta_{-1}}$ satisfies $f'_{\delta_{-1}}(t)=c$ with $f_{\delta_{-1}}(0)=0$. Hence, $f_{\delta_{-1}}(t)=ct$, and the variance in \eqref{def:gaussianlimit} simplifies to
\begin{align*}
\varfunc_{\delta_{-1}}(t) &= \sigma^2 \int_0^t \left(\frac{1 + c(T-t)}{1 + c(T-s)}\right)^2 ds = \sigma^2\frac{t(1 + c(T-t))}{1+cT}.
\end{align*}
Recall now from the proof of Lemma \ref{lemma_degree} that $\Var(\mu_{G_n})=1/\delta(G_n)$.
Letting $C_0=\tfrac12 c^2t^2+\tfrac16 c^3T^3$,
the bounds of Proposition \ref{pr:f-convergence} show that  $|f_{G_n}(t)-ct|\le C_0/\delta(G_n)$. The function $(1-\lambda x)/(1-\lambda y)$ is Lipschitz in $(x,y) \in [0,cT]^2$, uniformly in $\lambda \in [-2,0]$. From this it is straightforward to argue that $|\varfunc_{G_n}(t)-\varfunc_{\delta_{-1}}(t)|\le C_1/\delta(G_n)$ for some constant $C_1$, and similarly for the convergence of the value using the identity \eqref{def:valueidentity}.
Lastly,  the SDE \eqref{def:denseSDE} admits the solution $X(t)=\int_0^t\frac{1+c(T-t)}{1+c(T-s)}\sigma dW(s)$, and it is then straightforward to identify the Gaussian law $X(t)\sim\mathcal{N}(0,\varfunc_{\delta_{-1}}(t))$.

\subsubsection{Cycle graph case: Proof of Proposition \ref{prop:cyclegraph}.}\label{se:cyclegraphproof}

We begin by simplifying the expression of $Q'$, with $Q=Q_\mu$ defined as in \ref{def:Qcycle}. Differentiate under the integral sign to get
\begin{equation*}
Q'(x) = Q(x) \int_0^1 \frac{1 - \cos{\left(2\pi u\right)}}{1 + x - x \cos{\left(2\pi u\right)}}du =  \frac{Q(x)}{2\pi} \int_{-\pi}^{\pi} \frac{1 - \cos{\left(u\right)}}{1 + x - x \cos{\left(u\right)}}du.
\end{equation*}
We will first assume that $x > 0$ and then show that the formula is still valid for $x =0$. We perform the change of variable $t=\tan(u/2)$, using $\cos(u) = \frac{1-t^2}{1+t^2}$, to get
\begin{align*}
\int_{-\pi}^{\pi} \frac{1 - \cos(u)}{1 + x - x \cos(u)}du &= \int_{-\infty}^{\infty} \frac{1 - \frac{1-t^2}{1+t^2}}{1 + x - x \frac{1-t^2}{1+t^2}} \frac{2}{1+ t^2} dt =  \int_{-\infty}^{\infty} \frac{4t^2}{(1 + t^2(1+2x))(1+t^2)}dt \\
	&= \frac{2}{x}\int_{-\infty}^{\infty} \left(\frac{1}{1+t^2} - \frac{1}{1+ t^2(1+ 2x)} \right)dt\\
	&= \frac{2}{x} \left[\arctan(t) - \frac{1}{\sqrt{1+2x}}\arctan(t\sqrt{1+2x}) \right]_{t=-\infty}^{\infty} \\
	&= 2\pi \frac{1 + 2x - \sqrt{1+2x}}{x(1+2x)} =: 2\pi h(x).
\end{align*}
We thus find $Q'(x)=Q(x)h(x)$ for $x > 0$, and if we define $h(0):=1$ then the formula extends by continuity to $x=0$. Using $Q(0)=1$, we find $Q(x)=\exp\int_0^x h(u)du$, and we compute this integral using the change of variables $v=\sqrt{1+2u}$:
\begin{align*}
\int_0^x h(u) du &= \int_0^x \frac{1 + 2u - \sqrt{1+2u}}{u(1+2u)}du = \int_1^{\sqrt{1+2x}} \frac{2}{v+1} dv = \log\left(\frac{1}{2}(\sqrt{1+2x} + 1 + x) \right).
\end{align*}
Therefore
\begin{equation*}
Q(x) = \frac{1}{2}(\sqrt{1+2x} + 1 + x), \quad \text{and} \quad Q'(x) = \frac{1}{2}\left(\frac{1}{\sqrt{1+2x}} + 1\right).
\end{equation*}
Recall from Proposition \ref{prop:cyclegraph} that we defined $\Phi(x)=\log(1+\sqrt{1+2x}) -\sqrt{1+2x} + x + \tfrac12$, and $f(t):=f_\mu(t):=\Phi^{-1}\left(\log 2 +\frac{ct-1}{2} \right)$. It remains to show that $f$ satisfies the desired ODE $f'(t)=cQ'(f(t))$ with $f(0)=0$, where we recall that this ODE is well-posed by Proposition \ref{pr:f-convergence}.
Note that $\Phi$ is continuous and increasing on $\R_+$ and maps $\R_+$ to $[\log 2 - 1/2, \infty)$. Hence, the inverse $\Phi^{-1}$ is well defined on $[\log 2 - 1/2, \infty)$ with $f(0)=\Phi^{-1}(\log 2 - 1/2)=0$. Straightforward calculus yields $\Phi'(x)=\frac{\sqrt{1+2x}}{1+\sqrt{1+2x}}$ and thus $Q'(x)=\frac{1}{2\Phi'(x)}$, and we find that indeed $f'(t) = \frac{c}{2\Phi'\left(f(t)\right)} = cQ'(f(t))$,
which concludes the proof of Proposition \ref{prop:cyclegraph}.

\bibliographystyle{amsplain}
\bibliography{biblio_network}

\end{document}